%&AMSTeX
\input amstex
\documentstyle {amsppt}
\pagewidth{12.5 cm}\pageheight{19 cm}\magnification\magstep1
\topmatter
\title The $G$-stable pieces of the wonderful compactification \endtitle
\author Xuhua He\endauthor
\address Department of Mathematics, M.I.T., Cambridge, MA 02139\endaddress
\email xuhua\@mit.edu, hugo\@math.mit.edu \endemail
 \subjclass 20G15, 14L30 \endsubjclass

\abstract Let $G$ be a connected, simple algebraic group over an
algebraically closed field. There is a partition of the wonderful
compactification $\bar{G}$ of $G$ into finite many $G$-stable
pieces, which was introduced by Lusztig. In this paper, we will
investigate the closure of any $G$-stable piece in $\bar{G}$. We
will show that the closure is a disjoint union of some $G$-stable
pieces, which was first conjectured by Lusztig. We will also prove
the existence of cellular decomposition if the closure contains
finitely many $G$-orbits.
\endabstract

\endtopmatter
\document

\define\po{\text{\rm pos}}

\define\supp{\text{\rm supp}}

\define\Lie{\text{\rm Lie}}

\define\Ad{\text{\rm Ad}}
\redefine\i{^{-1}}
\redefine\ge{\geqslant}
\redefine\le{\leqslant}

\define\cb{\Cal B}

\define\ci{\Cal I}

\define\cp{\Cal P}

\define\cv{\Cal V}

\define\a{\alpha}
\redefine\b{\beta}

\define\g{\gamma}

\redefine\d{\delta}
\redefine\D{\Delta}

\define\p{\pi}

\redefine\v{\vartheta}

\define\tz{\tilde Z}

\head Introduction \endhead

An adjoint semi-simple group $G$ has a ``wonderful''
compactification $\bar{G}$, introduced by De Concini and Procesi
in \cite{DP}. The variety $\bar{G}$ is a smooth variety with $G
\times G$ action. Denote by $G_{diag}$, the image of the diagonal
embedding of $G$ into $G \times G$. The $G_{diag}$-orbits of
$\bar{G}$ were studied by Lusztig in \cite{L4}. He introduced a
partition of $\bar{G}$ into finitely many $G$-stable pieces. The
$G$-orbits on each piece can be described explicitly. Based on the
partition, he established the theory of ``parabolic character
sheaves'' on $\bar{G}$.

The main results of this paper concern the closure of the
$G$-stable pieces. The closure of each piece is a union of some
other pieces and if the closure contains finitely many $G$-orbits,
then it admits a cellular decomposition. I believe that our
results are necessary ingredients for establishing the
(conjectural) Kazhdan-Lusztig theory on the ``Parabolic Character
Sheaves'' on $\bar{G}$.

We now review the content of this paper in more detail.

In section 1, we recall the definition of $G$-stable pieces in
\cite{L4} and establish some basic results. The pieces are indexed
by the pairs $\ci=\{(J, w)\}$, where $J$ is a subset of the simple
roots and $w$ is an element of the Weyl group $W$, which has
minimal length in the coset $w W_J$. One interesting result is
that any $G$-stable piece is the minimal $G$-stable subset that
contains a particular $B \times B$-orbit, where $B$ is the Borel
subgroup. The closure of any $B \times B$-orbit in $\bar{G}$ was
studied by Springer in \cite{S}. Based on his result and the
relations between $G$-stable pieces and $B \times B$-orbits, we
are able to investigate the closure of the $G$-stable pieces.

In section 2, we recall the definition of the ``wonderful''
compactification and introduce ``compactification through the
fibres'', a technique tool that will be used to prove the
existence of cellular decomposition. In section 3, we describe a
partial order on $\ci$, which is the partial order that
corresponds to the closure relation of the $G$-stable piece, as we
will see in section 4. In section 4, we also discuss the closure
of any $G$-stable piece that appears in \cite{L3}.

In section 5, we discuss the existence of cellular decomposition.
Each piece does not have a cellular decomposition. However, a
union of certain pieces has a cellular decomposition. (This is
motivated by Springer in \cite{S}, in which he showed that a union
of certain $B \times B$-orbits is isomorphic to an affine space.)
In fact, if the closure contains finitely many $G$-orbits, then it
has a cellular decomposition.

The methods work for arbitrary connected component of a
disconnected algebraic group with identity component $G$. The
results for that component is just a ``twisted'' version of the
results for $G$ itself.

\head 1. The $G$-stable pieces \endhead

\subhead 1.1 \endsubhead In the sequel $G$ is a connected,
semi-simple algebraic group of adjoint type over an algebraically
closed field. Let $B$ be a Borel subgroup of $G$, $B^-$ be the
opposite Borel subgroup and $T=B \cap B^-$. Let $(\a_i)_{i \in I}$
be the set of simple roots. For $i \in I$, we denote by $s_i$ the
corresponding simple reflection. For any element $w$ in the Weyl
group $W=N(T)/T$, we will choose a representative $\dot w$ in
$N(T)$ in the same way as in \cite{L1, 1.1}.

For $J \subset I$, let $P_J \supset B$ be the standard parabolic
subgroup defined by $J$ and $P^-_J \supset B^-$ be the opposite of
$P_J$. Set $L_J=P_J \cap P^-_J$. Then $L_J$ is a Levi subgroup of
$P_J$ and $P^-_J$. We denote by $\Phi_J$ the set of roots that are
linear combination of $\{(a_j)_{j \in J}\}$. Let $Z_J$ be the
center of $L_J$ and $G_J=L_J/Z_J$ be its adjoint group. We denote
by $\p_{P_J}$ (resp. $\p_{P^-_J}$) the projection of $P_J$ (resp.
$P^-_J$) onto $G_J$.

For any $J \subset I$, let $\cp_J$ be the set of parabolic
subgroups conjugate to $P_J$. We will write $\cb$ for
$\cp_{\varnothing}$.

For any subset $J$ of $I$, let $W_J$ be the subgroup of $W$
generated by $\{s_j \mid j \in J\}$ and $W^J$ (resp. $^J W$) be
the set of minimal length coset representatives of $W/W_J$ (resp.
$W_J \backslash W$). Let $w^J_0$ be the unique element of maximal
length in $W_J$. (We will simply write $w^I_0$ as $w_0$.) For $J,
K \subset I$, we write $^J W^K$ for $^J W \cap W^K$.

For $w \in W$, we denote by $\supp(w) \subset I$ the set of simple
roots whose associated simple reflections occur in some (or
equivalently, any) reduced decomposition of $w$.

For $J, K \subset I$, $P \in \cp_J, Q \in \cp_K$ and $u \in {}^J
W^K$, we write $\po(P, Q)=u$ if there exists $g \in G$, such that
$^g P=P_J$ and $^g Q={}^{\dot u} P_K$.

For any parabolic subgroup $P$, we denote by $U_P$ its unipotent
radical. We will simply write $U$ for $U_B$ and $U^-$ for
$U_{B^-}$. For $J \subset I$, set $U_J=U \cap L_J$ and $U^-_J=U^-
\cap L_J$.

For any closed subgroup $H$ of $G$, we denote by $H_{diag}$ the
image of the diagonal embedding of $H$ in $G \times G$. For any
subgroup $H$ and $g \in G$, we write $^g H$ for $g H g\i$. For any
finite set $A$, we write $|A|$ for the cardinal of $X$.

\subhead 1.2 \endsubhead Let $\hat{G}$ be a possibly disconnected
reductive algebraic group over and algebraically closed field with
identity component $G$. Let $G^1$ be a fixed connected component
of $\hat{G}$. There exists an isomorphism $\d: W@>>>W$ such that
$\d(I)=I$ and $^g P \in \cp_{\d(J)}$ for $g \in G^1$ and $P \in
\cp_J$. There also exists $g_0 \in G^1$ such that $g_0$ normalizes
$T$ and $B$. Moreover, $g_0$ can be chosen in such a way that
$^{g_0} L_J=L_{\d(J)}$ for $J \subset I$. We will fix such $g_0$
in the rest of this paper.

In particular, if $G^1=G$, then $\d=id$, where $id$ is the
identity map. In this case, we choose $g_0$ to be the unit element
$1$ of $G$.

\subhead 1.3 \endsubhead We will follow the set-up of \cite{L4}.

Let $J, J' \subset I$ and $y \in {}^{J'} W^J$ be such that $\Ad(y)
\d(J)=J'$. For $P \in \cp_J$, $P' \in \cp_{J'}$, define $A_y(P,
P')=\{g \in G^1 \mid \po(P', {}^g P)=y\}$. Define $$\tz_{J, y,
\d}=\{(P, P', \g) \mid P \in \cp_J, P' \in \cp_{J'}, \g \in U_{P'}
\backslash A_y(P, P') /U_P\}$$ with $G \times G$ action defined by
$(g_1, g_2) (P, P', \g)=(^{g_2} P, {}^{g_1} P', g_1 \g g_2 \i)$.

By \cite{L4, 8.9}, $A_y(P, P')$ is a single $P', P$ double coset.
Thus $G \times G$ acts transitively on $Z_{J, y, \d}$.

Let $z=(P, P', \g) \in \tz_{J, y, \d}$. Then there exists $g \in
\g$ such that $^g P$ contains some Levi of $P \cap P'$. Now set
$P_1=g \i (^g P)^{({P'}^P)} g$, $P'_1=P^{P'}$. Define $$\a(P, P',
\g)=(P_1, P'_1, U_{P'_1} g U_{P_1}).$$ By \cite{L4, 8.11}, The map
$\a$ doesn't depend on the choice of $g$.

To $z=(P,P',\g)\in \tz_{J, y, \d}$, we associate a sequence $(J_k,
J'_k, u_k, y_k, P_k, P'_k, \g_k)_{k \ge 0}$ with $J_k, J'_k
\subset I$, $u_k \in W$, $y_k \in
{}^{J'_k}W^{\d(J_k)},\Ad(y_k)\d(J_k)=J'_k$, $P_k \in \cp_{J_k},
P'_k \in \cp_{J'_k}, \g_k=U_{P'_k} g U_{P_k}$ for some $g \in G$
satisfies $\po(P'_k, {}^g P_k)=y_k$. The sequence is defined as
follows.
$$P_0=P, P'_0=P', \g_0=\g, J_0=J, J'_0=J', u_0=\po(P'_0,P_0), y_0=y.$$
Assume that $k\ge 1$, that $P_m, P'_m, \g_m, J_m, J'_m, u_m, y_m$
are already defined for $m<k$ and that $u_m=\po(P'_m,P_m), P_m \in
\cp_{J_m}, P'_m \in \cp_{J'_m}$ for $m<k$. Let
$$J_k=J_{k-1} \cap \d \i \Ad(y_{k-1}\i u_{k-1})J_{k-1},
J'_k=J_{k-1} \cap \Ad(u_{k-1}\i y_{k-1}) \d(J_{k-1}),$$
$$(P_k, P'_k, \g_k)=\a(P_{k-1}, P'_{k-1}, \g_{k-1}) \in \tz_{J_k, y_k, \d} (\hbox{see \cite {L4, 8.10}}),$$
$$u_k=\po(P'_k,P_k),y_k=u_{k-1}\i y_{k-1},\g_k=U_{P'_k}g_{k-1}U_{P_k}.$$

It is known that the sequence is well defined. Moreover, for
sufficient large $n$, we have that
$J_n=J'_n=J_{n+1}=J'_{n+1}=\cdots=J_{\infty}$,
$u_n=u_{n+1}=\cdots=1$, $y_n=y_{n+1}=\cdots=y_{\infty}$,
$P_n=P_{n+1}=\cdots=P_{\infty}$,
$P'_n=P'_{n+1}=\cdots=P'_{\infty}$ and
$\g_n=\g_{n+1}=\cdots=\g_{\infty}$. Now we set $\b(z)=y_{\infty}$.
Then we have that $\b(z) \in W^{\d(J)}$. By \cite{L4, 8.18} and
\cite{L3, 2.5}, the sequence $(J_n, J'_n, u_n, y_n)_{n \ge 0}$ is
uniquely determined by $\b(z)$ and $y$.

For $w \in W^{\d(J)}$, set $$\tz^w_{J, y, \d}=\{z \in \tz_{J, y,
\d} \mid \b(z)=w\}.$$

Then $(\tz^w_{J, y, \d})_{w \in W^{\d(J)}}$ is a partition of
$\tz_{J, y, \d}$ into locally closed $G$-stable subvarieties. We
call $(\tz^w_{J, y, \d})_{w \in W^{\d(J)}}$ the $G$-stable pieces
of $\tz_{J, y, \d}$. For $w \in W^{\d(J)}$, let $(J_n, J'_n, u_n,
y_n)_{n \ge 0}$ be the sequence determined by $w$ and $y$. The
restriction of the map $\a$ on $\tz^w_{J, y, \d}$ is a
$G$-equivariant morphism from $\tz^w_{J, y, \d}$ onto $\tz^w_{J_1,
y_1, \d}$. We also denote this morphism by $\a$. It is known that
$\a$ induces a bijection from the set of $G$-orbits on $\tz^w_{J,
y, \d}$ to the set of $G$-orbits on $\tz^w_{J_1, y_1, \d}$.

We have a consequence $\tz^w_{J, y, \d}@>\a>>\tz^w_{J_1, y_1,
\d}@>\a>>\tz^w_{J_2, y_2, \d}@>\a>>\cdots$. For sufficiently large
$n$, $\v=\a^n: \tz^w_{J, y, \d}@>>>\tz^w_{J_{\infty}, w, \d}$ is
independent of the choice of $n$ and is a $G$-equivariant
morphism. Moreover, $\v$ induces a bijection from the set of
$G$-orbits on $\tz^w_{J, y, \d}$ to the set of $G$-orbits on
$\tz^w_{J_{\infty}, w, \d}$.

${}$

In the rest of this section, we will fix $J, y, \d, w$ and
$J_{\infty}$. First, we will give an explicit description of
$J_{\infty}$ in terms of $J, \d$ and $w$.

\proclaim{Lemma 1.4} Keep the notion of 1.3. Then
$$J_{\infty}=\max\{K \subset J \mid \Ad(w)\d(K)=K\}.$$
\endproclaim

Proof. Set $v=y_1 w \i$. By \cite{H, 2.2}, $v \in W_J$. Now $J_1=J
\cap \d \i \Ad(y_1 \i) J$. Thus $\Phi_{\d(J_1)} \subset \Ad(y_1
\i) \Phi_J=\Ad(w \i) \Ad(v \i) \Phi_J=\Ad(w \i) \Phi_J$.

Let $i \in J$. Assume that $\a_{\d(i)} \in \Ad(y_1 \i) \Phi_J$.
Then $\a_{\d(i)}=\Ad(y_1 \i) \a=\Ad(y \i) \Ad(u_0) \a$ for some
$\a \in \Phi_J$. Then $\a_{\Ad(y) \d(i)}=\Ad(u_0) \a$. Note that
$\a_{\Ad(y) \d(i)}$ is a simple root and $u_0 \in W^J$. Then
$\a=\a_j$ for some $j \in J$. Hence $i=\d \i \Ad(y_1 \i) j$.
Therefore, $i \in J \cap \d \i \Ad(y_1 \i) J=J_1$. So
$$J_1=\max\{K \subset J \mid \Phi_{\d(K)} \subset \Ad(w \i)
\Phi_J\}.$$

Set $J'_{\infty}=\max\{K \subset J \mid \Ad(w)\d(K)=K\}$. Then
$J'_{\infty} \subset J$. Moreover, $\Phi_{\d(J'_{\infty})}=\Ad(w
\i) \Phi_{J'_{\infty}} \subset \Ad(w \i) \Phi_J$. Thus
$J'_{\infty} \subset J_1$. We can show by induction that
$J'_{\infty} \subset J_n$ for all $n$. Thus $J'_{\infty} \subset
J_{\infty}$. By the definition, $J_{\infty}=J_{\infty} \cap \d \i
\Ad(w \i) J_{\infty}$. Thus $\Ad(w) \d(J_{\infty})=J_{\infty}$. So
$J_{\infty}=J'_{\infty}$. The lemma is proved. \qed

\subhead 1.5 \endsubhead Now set $\tilde{h}_{J, y, \d}=(P_J,
{}^{\dot y \i} P_{J'}, U_{^{\dot y \i} P_{J'}} g_0 U_{P_J}) \in
Z_{J, y, \d}$. For $w \in W^{\d(J)}$ and $v \in W$, set
$\widetilde{[J, w, v]}_{y, \d}=(B \times B) (\dot w, \dot v) \cdot
\tilde{h}_{J, y, \d}$. Then we have the following result.

\proclaim{Lemma 1.6} Keep the notion of 1.3. Let $g \in P_{J_1}$.
Set $z=(\dot w, g) \cdot \tilde{h}_{J, y, \d}$ and $z'=(\dot w, g)
\cdot \tilde{h}_{J_1, y_1, \d}$. Then $\a(z)=z'$.
\endproclaim

Proof. Set $P=P_J$, $P'={}^{\dot w \dot y \i} P_{J'}$, $g_1=\dot w
g_0 g$ and $v=y_1 w \i$. Then $v \in W_J$.

By the proof of \cite{H, 2.3}, $^{\dot v \i} L_{J'_1}$ is a Levi
factor of $P \cap P'$ and $P^{P'}={}^{\dot v \i} P_{J'_1}$,
$(P')^P={}^{\dot w \dot y \i} P_{\Ad(y) \d(J_1)}$. Moreover,
$$^{\dot v \i} L_{J'_1}={}^{\dot w \dot y_1 \i} L_{J'_1}={}^{\dot w}
L_{\d(J_1)}={}^{\dot w g_0} L_{J_1} \subset {}^{\dot w g_0}
P_J={}^{\dot w g_0 g} P_J.$$

So $^{g_1} P$ contains some Levi of $P \cap P'$. We have that
$$\eqalignno{g_1 \i (^{g_1} P)^{(^{\dot w \dot y \i} P_{\Ad(y_1) \d(J_1)})}
g_1 &={}^{g \i} P^{(^{g_0 \i \dot y \i} P_{\Ad(y) \d(J_1)})}={}^{g
\i} P_{J \cap \d \i \Ad(y \i) \Ad(y) \d(J_1)} \cr &={}^{g \i}
P_{J_1}=P_{J_1}.}$$

Thus $\a(z)=z'$. The lemma is proved. \qed

\proclaim{Proposition 1.7} We have that $$\tz^w_{J, y,
\d}=G_{diag} \cdot \widetilde{[J, w, 1]}_{y, \d}=G_{diag} \cdot
(P_J, {}^{\dot w \dot y \i} P_{J'}, U_{^{\dot w \dot y \i} P_{J'}}
\dot w g_0 (B \cap L_{J_{\infty}}) U_{P_J}).$$
\endproclaim

Proof. It is easy to see that $\tz^w_{J_{\infty}, w, \d}=G_{diag}
(\dot w, L_{J_{\infty}}) \cdot \tilde{h}_{J_{\infty}, w, \d}$.
Thus for any $b \in B$, $\a^n \bigl((\dot w, b) \cdot
\tilde{h}_{J, y, \d} \bigr) \in \tz^w_{J_{\infty}, w, \d}$ for
sufficiently large $n$. Therefore, $G_{diag} (\dot w, B) \cdot
\tilde{h}_{J, y, \d} \subset \tz^w_{J, y, \d}$.

Note that $\dot w g_0$ normalizes $(L_{J_{\infty}})$ and
$(L_{J_{\infty}}) \cap B$. Thus $$\dot w g_0 L_{J_{\infty}}=\{l
\dot w g_0 b l \i \mid l \in L_{J_{\infty}}, b \in L_{J_{\infty}}
\cap B\}.$$ Hence any element in $\tz^w_{J_{\infty}, w, \d}$ is
$G$-conjugate to $(\dot w, l) \cdot \tilde{h}_{J_{\infty}, w, \d}$
for some $l \in L_{J_{\infty}} \cap B$. Now let $z \in \tz^w_{J,
y, \d}$. Then $\v(z)$ is $G$-conjugate to $(\dot w, l) \cdot
\tilde{h}_{J_{\infty}, w, \d}$ for some $l \in L_{J_{\infty}} \cap
B$. Set $z'=(\dot w, l) \cdot \tilde{h}_{J, y, \d} \in \tz^w_{J,
y, \d}$. Then $\v(z')$ lies in the same $G$-orbit as $\v(z)$.
Since $\v$ induces a bijection from the set of $G$-orbits on
$\tz^w_{J, y, \d}$ to the set of $G$-orbits on $\tz^w_{J_{\infty},
w, \d}$. Thus $z$ is $G$-conjugate to $z'$. So $\tz^w_{J, y,
\d}=G_{diag} (\dot w, B \cap L_{J_{\infty}})) \cdot \tilde{h}_{J,
y, \d}$. The proposition is proved. \qed

\subhead 1.8\endsubhead In \cite{L4, 8.20}, Lusztig showed that
$\tz^w_{J, y, \d}$ is an iterated affine space bundle over a fibre
bundle over $\cp_{J_{\infty}}$ with fibres isomorphic to
$L_{J_{\infty}}$. In 1.10, we will prove a similar (but more
explicit) result, which will be used to establish the cellular
decomposition. Before doing that, we include the following result
(see \cite{SL, page 26, lemma 4}) which we will use in the proof
of proposition 1.10 and proposition 5.5.

\proclaim{Lemma 1.9} Let $H$ be a closed subgroup of $G$ and
$\Phi: X \rightarrow G/H$ be a $G$-equivariant morphism from the
$G$-variety $X$ to the homogeneous space $G/H$. Let $E \subset X$
be the fiber $\Phi \i(H)$. Then $E$ will be stabilized by $H$ and
the map $\Psi: G \times_H E \rightarrow X$ sending $(g, e)$ to $g
\cdot e$ defines an isomorphism of $G$-varieties.
\endproclaim

\proclaim{Proposition 1.10} For $a \in W$, set $U_a=U \cap
{}^{\dot a} U^-$. Set $$\tilde{L}^w_{J, y, d}=(L_{J_{\infty}},
L_{J_{\infty}}) (\dot w, 1) \cdot \tilde{h}_{J, y, \d}.$$

Then we have the following results.

(1) $\tz^w_{J, y, \d}$ is isomorphic to $G \times_{P_{J_{\infty}}}
\bigl( (P_{J_{\infty}}) \cdot \tilde{L}^w_{J, y, \d} \bigr)$.

(2) $(P_{J_{\infty}}) \cdot \tilde{L}^w_{J, y, \d}=(B \times B)
\cdot \tilde{L}^w_{J, y, \d} \cong (U \cap {}^{\dot
w^{J_{\infty}}_0 \dot w \dot y \i \dot w^{J'}_0} U^-) \times
\tilde{L}^w_{J, y, \d}$, where $\tilde{L}^w_{J, y, \d}$ is
isomorphic to $L_{J_{\infty}}$.

(3) $G_{diag} (\dot w T, 1) \cdot \tilde{h}_{J, y, \d}$ is dense
in $\tz^w_{J, y, \d}$.
\endproclaim

Proof. It is easy to see that $\tilde{L}^w_{J, y, \d}=(\dot w,
L_{J_{\infty}}) \cdot \tilde{h}_{J, y, \d}$ is isomorphic to
$L_{J_{\infty}}$. By 1.7, $\tz^w_{J, y, \d}=G_{diag} \cdot
\tilde{L}^w_{J, y, \d}$. Consider the $G$-equivariant map $p:
\tz^w_{J_{\infty}, w, \d}@>>> \cp_{J_{\infty}}$ defined by $p(P,
P, \g)=P$ for $(P, P, \g) \in \tz^w_{J_{\infty}, w, \d}$. For $l
\in L_{J_{\infty}}$ and $g \in G$, if $p \circ \v \bigl( (g, g)
(\dot w, l) \cdot \tilde{h}_{J, y, \d} \bigr)=P_{J_{\infty}}$,
then $g \in P_{J_{\infty}}$. Thus $(P_{J_{\infty}})_{diag} \cdot
\tilde{L}^w_{J, y, \d}=(p \circ \v) \i (P_{J_{\infty}})$.

Assume that $(g, g) (\dot w, l_1) \cdot \tilde{h}^w_{J, y,
\d}=(\dot w, l_2) \cdot \tilde{h}^w_{J, y, \d}$ for $g \in G$ and
$l_1, l_2 \in L_{J_{\infty}}$. Then $^g P_{J_{\infty}}=p \circ \v
\bigl((g, g) (\dot w, l_1) \cdot \tilde{h}^w_{J, y, \d} \bigr)=p
\circ \v \bigl((\dot w, l_2) \cdot \tilde{h}^w_{J, y, \d}
\bigr)=P_{J_{\infty}}$. So $g \in P_{J_{\infty}}$. Part (1) is
proved.

We have that $(B \times B) \cdot \tilde{L}^w_{J, y, \d}=(B)_{diag}
(1, B) \cdot \tilde{L}^w_{J, y, \d}$ and $p \circ \v \bigl( (1, B)
\cdot \tilde{L}^w_{J, y, \d} \bigr)=P_{J_{\infty}}$. Thus $(B
\times B) \cdot \tilde{L}^w_{J, y, \d} \subset (B)_{diag}
(P_{J_{\infty}})_{diag} \tilde{L}^w_{J, y,
\d}=(P_{J_{\infty}})_{diag} \tilde{L}^w_{J, y, \d}$. On the other
hand, $(P_{J_{\infty}})_{diag} \tilde{L}^w_{J, y, \d} \subset
(P_{J_{\infty}}, P_{J_{\infty}}) \cdot \tilde{L}^w_{J, y, \d}=(B
\times B) \cdot \tilde{L}^w_{J, y, \d}$. Hence $(P_{J_{\infty}})
\cdot \tilde{L}^w_{J, y, \d}=(B \times B) \cdot \tilde{L}^w_{J, y,
\d}$.

Now consider $\pi: (U \cap {}^{\dot w^{J_{\infty}}_0 \dot w \dot y
\i \dot w^{J'}_0} U^-) \times \tilde{L}^w_{J, y, \d}@>>>(B \times
B) \cdot \tilde{L}^w_{J, y, \d}$ defined by $\pi(u, l)=(u, 1) l$
for $u \in U \cap {}^{\dot w^{J_{\infty}}_0 \dot w \dot y \i \dot
w^{J'}_0} U^-$ and $l \in \tilde{L}^w_{J, y, \d}$.

Note that $(1, B L_{J_{\infty}}) \cdot \tilde{h}_{J, y, \d}=(1,
U_{P_J} L_{J_{\infty}} U_J) \cdot \tilde{h}_{J, y, \d}=(U_{\d(J)},
L_{J_{\infty}}) \cdot \tilde{h}_{J, y, \d}$. Since $w \in
W^{\d(J)}$, $B \dot w U_{\d(J)}=B \dot w=U_{P_{J_{\infty}}}
L_{J_{\infty}} \dot w$. Hence $$(B \dot w, B L_{J_{\infty}}) \cdot
\tilde{h}_{J, y, \d}=(U_{P_{J_{\infty}}} \dot w
L_{\d(J_{\infty})}, L_{J_{\infty}}) \cdot \tilde{h}_{J, y,
\d}=(U_{P_{J_{\infty}}} \dot w, L_{J_{\infty}}) \cdot
\tilde{h}_{J, y, \d}.$$ Since $w y \i \in W^{J'}$ and $\Ad(y w \i)
J_{\infty} \subset J'$, then
$$\eqalignno{U_{P_{J_{\infty}}} &=(U_{P_{J_{\infty}}} \cap {}^{\dot w \dot y \i
\dot w^{J'}_0} U^-) (U_{P_{J_{\infty}}} \cap {}^{\dot w \dot y \i}
U_{P_{J'}}) \cr &=(U \cap {}^{\dot w^{J_{\infty}}_0 \dot w \dot y
\i \dot w^{J'}_0} U^-) (U_{P_{J_{\infty}}} \cap {}^{\dot w \dot y
\i} U_{P_{J'}}).}$$

Therefore, $(B \dot w, B L_{J_{\infty}}) \cdot \tilde{h}_{J, y,
\d}=\bigl((U \cap {}^{\dot w^{J_{\infty}}_0 \dot w \dot y \i \dot
w^{J'}_0} U^-) \dot w, L_{J_{\infty}} \bigr) \cdot \tilde{h}_{J,
y, \d}$. So $\pi$ is surjective.

Let $u\in U \cap {}^{\dot w^{J_{\infty}}_0 \dot w \dot y \i \dot
w^{J'}_0} U^-$ and $l_1, l_2 \in L_{J_{\infty}}$. Assume that $(u,
1) (\dot w, l_1) \cdot \tilde{h}_{J, y, \d}=(\dot w, l_2) \cdot
\tilde{h}_{J, y, \d}$. Note that the isotropy subgroup of $G
\times G$ at the point $(\dot w, 1) \cdot \tilde{h}_{J, y, \d}$ is
$\{(U_{^{\dot w \dot y \i} P_{J'}} l', U_{P_J} g_0 \i \dot w \i l'
\dot w g_0) \mid l' \in {}^{\dot w} L_{\d(J)}\}$.

Thus $u \in U_{^{\dot w \dot y \i} P_{J'}} l'$ and $l_2 \i l_1 \in
U_{P_J} g_0 \i \dot w \i l' \dot w g_0$ for some $l' \in {}^{\dot
w} L_{\d(J)}$. Then $l' \in L_{J_{\infty}}$ and $u=1$. Thus $\pi$
is injective.

In fact, we can show that the bijective morphism $\pi$ is an
isomorphism. The verification is omitted.

Part (3) can be proved in the same way as in \cite{H, 2.7}. \qed

\subhead 1.11\endsubhead For $P\in\cp_J$, let $H_P$ be the inverse
image of the connected center of $P/U_P$ under $P@>>>P/U_P$. We
can regard $H_P/U_P$ as a single torus $\D_J$ independent of $P$.
Then $\D_J$ acts (freely) in the natural way on $\tz_{J, y, \d}$
and the action commutes with the action of $G$. Moreover, each
piece $\tz^w_{J, y, \d}$ is $\D_J$-stable.

Define
$$\eqalignno{Z_{J, y, \d} &=\{(P, P', \g) \mid P \in \cp_J, P' \in \cp_{J'},
\g \in H_{P'} \backslash A_y(P, P') /U_P\} \cr &=\{(P, P', \g)
\mid P \in \cp_J, P' \in \cp_{J'}, \g \in U_{P'} \backslash A_y(P,
P') /H_P\}}$$ with $G \times G$ action defined by $(g_1, g_2) (P,
P', \g)=(^{g_2} P, {}^{g_1} P', g_1 \g g_2 \i)$.

Then $Z_{J, y, \d}$ and $\D_J \backslash \tz_{J, y, \d}$ can be
identified in the natural way as varieties with $G$-action. Set
$Z^w_{J, y, \d}=\D_J \backslash \tz^w_{J, y, \d}$. Then $$Z_{J, y,
\d}=\sqcup_{w \in W^{\d(J)}} Z^w_{J, y, \d}.$$

We call $(Z^w_{J, y, \d})_{w \in W^{\d(J)}}$ the $G$-stable pieces
of $Z_{J, y, \d}$. Set $$\eqalignno{h_{J, y, \d} &=(P_J, {}^{\dot
y \i} P_{J'}, H_{^{\dot y \i} P_{J'}} g_0 U_{P_J}) \in Z_{J, y,
\d}, \cr L^w_{J, y, \d} &=(L_{J_{\infty}}, L_{J_{\infty}}) (\dot
w, 1) \cdot h_{J, y, \d}.}$$ For $w \in W^{\d(J)}$ and $v \in W$,
set $[J, w, v]_{y, \d}=(B \times B) (\dot w, \dot v) \cdot h_{J,
y, \d}$. Then as a consequence of 1.7 and 1.10, we have the
following result.

\proclaim{Proposition 1.12} For $w \in W^{\d(J)}$, we have that

(1) $Z^w_{J, y, \d}=G_{diag} \cdot [J, w, 1]_{y, \d}$.

(2) $Z^w_{J, y, \d}$ is isomorphic to $G \times_{P_{J_{\infty}}}
\bigl( (P_{J_{\infty}}) \cdot L^w_{J, y, \d} \bigr)$.

(3) $(P_{J_{\infty}}) \cdot L^w_{J, y, \d}=(B \times B) \cdot
L^w_{J, y, \d} \cong (U \cap {}^{\dot w^{J_{\infty}}_0 \dot w \dot
y \i \dot w^{J'}_0} U^-) \times L^w_{J, y, \d}$, where $L^w_{J, y,
\d}$ is isomorphic to $L_{J_{\infty}} / Z(L_J)$.

(4) $G_{diag} (\dot w T, 1) \cdot h_{J, y, \d}$ is dense in
$Z^w_{J, y, \d}$.
\endproclaim

\head 2. Compactification through the fibres \endhead

\subhead 2.1\endsubhead For any connected, semi-simple algebraic
group of adjoint type, De Concini and Procesi introduced its
wonderful compactification $\bar{G}$(see \cite{DP}). It is an
irreducible, projective smooth $G \times G$-variety. The $G \times
G$-orbits $Z_J$ of $\bar{G}$ are indexed by the subsets $J$ of
$I$. Moreover, $Z_J=(G \times G) \times_{P^-_J \times P_J} G_J$,
where $P^-_J \times P_J$ acts on the right on $G \times G$ and on
the left on $G_J$ by $(q, p) \cdot z=\p_{P^-_J}(q) z \p_{P_J}(p)$.
Denote by $h_J$ the image of $(1, 1, 1)$ in $Z_J$. We will
identify $Z_J$ with $Z_{J, w_0 w^J_0, id}$ and $h_J$ with $h_{J,
w_0 w^J_0, id}$, where $id$ is the identity map on $I$ (see
\cite{H, 2.5}).

Let us consider the $B \times B$-orbits on $\bar{G}$. For any $J
\subset I$, $u \in W^J$ and $v \in W$, set $[J, u, v]=(B \times B)
(\dot{u}, \dot{v}) \cdot h_J$. Then $\bar{G}=\bigsqcup\limits_{J
\subset I} \bigsqcup\limits_{x \in W^J, w \in W} [J, x, w]$. The
following result is due to Springer (see \cite{S, 2.4}).

\proclaim{Theorem} Let $x \in W^J$, $x' \in W^K$, $w, w' \in W$.
Then $[K, x', w']$ is contained in the closure of $[J, x, w]$ if
and only if $K \subset J$ and there exists $u \in W_K, v \in W_J
\cap W^K$ with $x v u \i \le x'$, $w' u \le w v$ and $l(w
v)=l(w)+l(v)$. In particular, for $J \subset I$ and $w \in W^J$,
the closure of $[J, w, 1]$ in $\bar{G}$ is $\sqcup_{K \subset J}
\sqcup_{x \in W^K, u \in W_J, \hbox{ and } x \ge w u} [K, x, u]$.
\endproclaim

\subhead 2.2\endsubhead We have defined $Z_{J, y, \d}$ in 1.11. As
we have seen, $Z_{J, y, \d}$ is a locally trivial fibre bundle
over $\cp_J \times \cp_{J'}$ with fibres isomorphic to
$L_J/Z(L_J)$. Note that $L_J/Z(L_J)$ is a connected, semi-simple
algebraic group of adjoint type. Thus we can define the wonderful
compactification $\overline{L_J/Z(L_J)}$ of $L_J/Z(L_J)$. In this
section, we will define $\underline{Z_{J, y, \d}}$, which is a
locally trivial fibre bundle over $\cp_J \times \cp_{J'}$ with
fibres isomorphic to $\overline{L_J/Z(L_J)}$.

\subhead 2.3\endsubhead We keep the notation of 1.3. Fix $g \in
A_y(P, P')$. Then $A_y(P, P') g \i=P' U_{^g P}$ (see \cite{L4,
8.9}). Set $$L_{P, P', g}={}^g P \cap P'/H_{^g P \cap P'}.$$

Let $\Phi_g: H_{P'} \backslash A_y(P, P') /H_P@>>>L_{P, P', g}$ be
the morphism defined by $$H_{P'} \backslash A_y(P, P') /H_P@>\cdot
g \i>>H_{P'} \backslash A_y(P, P') g \i /H_{^g P}@<i<<L_{P, P',
g},$$ where $i$ is the obvious isomorphism.

The $P \times P'$ action on $H_{P'} \backslash A_y(P, P') /H_P$
induces a $P \times P'$ action on $L_{P, P', g}$. Now for $g, g'
\in A_y(P, P')$, set $\Phi_{g, g'}=\Phi_{g'} \Phi_g \i: L_{P, P',
g}@>\simeq>>L_{P, P', g'}$. Then $\Phi_{g, g'}$ is compatible with
the $P \times P'$ action. Moreover, $(L_{P, P', g}, \Phi_{g, g'})$
forms an inverse system and $$H_{P'} \backslash A_y(P, P')
/H_P=\lim\limits_{\leftarrow} L_{P, P', g}.$$

Note that $L_{P, P', g}$ is a semi-simple group of adjoint type.
Then we can define the De Concini-Procesi compactification
$\overline{L_{P, P', g}}$ of $L_{P, P', g}$. The $P \times P'$
action on $L_{P, P', g}$ can be extended in the unique way to a $P
\times P'$ action on $\overline{L_{P, P', g}}$. The isomorphism
$\Phi_{g, g'}: L_{P, P', g}@>\simeq>>L_{P, P', g'}$ can be
extended in the unique way to an isomorphism from $\overline{L_{P,
P', g}}$ onto $\overline{L_{P, P', g'}}$. We will also denote this
isomorphism by $\Phi_{g, g'}$. It is easy to see that this
isomorphism is compatible with the $P \times P'$ action. Now
$\overline{(L_{P, P', g}}, \Phi_{g, g'})$ forms an inverse system.
Define
$$\overline{H_{P'} \backslash A_y(P, P') /H_P}=\lim\limits_{\leftarrow} \overline{L_{P,
P', g}}.$$

We also obtain a $P \times P'$ action on $\overline{H_{P'}
\backslash A_y(P, P') /H_P}$. Thus we can identify
$\overline{H_{P'} \backslash A_y(P, P') /H_P}$ with
$\overline{L_{P, P', g}} g$ as varieties with $P \times P'$
action.

\subhead Remark \endsubhead $\overline{H_{P'} \backslash A_y(P,
P') /H_P}$ is isomorphic to $\overline{L_{P, P', g}}$ as a
variety. However, we are also concerned with the $P' \times P$
action. In this case, $\overline{H_{P'} \backslash A_y(P, P')
/H_P}$ is regarded as $\overline{L_{P, P', g}} g$ with ``twisted''
$P' \times P$ action.

\subhead 2.4\endsubhead In this section, we will consider a
special case, namely, $P=P'=G^0$. In this case, $A_y(P, P')=G^1$
and we will identify $H_G \backslash A_y(G, G) /H_G$ with $G^1$.

Let $\cv_G$ be the projective variety whose points are the
$\dim(G)$-dimensional Lie subalgebras of $\Lie(G \times G)$. The
$\hat{G} \times \hat{G}$ action on $\Lie(G \times G)$ which is
defined by $(g_1, g_2) \cdot (a, b)=(\Ad(g_2) a, \Ad(g_1) b)$ for
$g_1, g_2 \in \hat{G}$ and $a, b \in \Lie(G)$ induces a $\hat{G}
\times \hat{G}$ action on $\cv_G$. To each $g \in \hat{G}$, we
associate a $\dim(G)$-dimensional subspace $V_g=\{(a, \Ad(g) a)
\mid a \in \Lie(G)\}$ of $\Lie(G \times G)$. Then $V_{g_1 g g_2
\i}=(g_1, g_2) \cdot V_g$ for $g_1, g, g_2 \in \hat{G}$ and $g
\mapsto V_g$ is an embedding $G^1 \subset \cv_G$. We denote the
image by $i(G^1)$.

If $G^1=G$, then the closure of $i(G)$ in $\cv_G$ is $\bar{G}$
(see \cite{DP}). Note that $V_{g g_0}=(1, g_0 \i) V_g$ for all $g
\in G$. Thus $i(G^1)=(1, g_0 \i) i(G)$. Hence the closure of
$i(G^1)$ in $\cv_G$ is $(1, g_0 \i) \bar{G}$, which is just
$\bar{G^1}$ defined above.

\subhead Remark \endsubhead In \cite{L4, 12.3}, Lusztig defined
the compactification of $G^1$ to be the closure of $i(G^1)$ in
$\cv_G$. As we have seen, our definition coincides with his
definition.

\subhead 2.5\endsubhead In \cite{L4, 12.3}, Lusztig showed that
$$\bar{G^1}=\sqcup_{J \subset I} Z_{J, w_0 w^{\d(J)}_0, \d}=\sqcup_{J \subset I}
\sqcup_{w \in W^{\d(J)}} Z^w_{J, w_0 w^{\d(J)}_0, \d},$$ where
the base point $h_{J, w_0 w^{\d(J)}_0, \d}=(P_J, P^-_{\d(J)},
H_{P^-_{\d(J)}} g_0 H_{P_J})$ is identified with the
$\dim(G)$-dimensional subalgebra $\{(l u, g_0 l g_0 \i u') \mid l
\in L_J, u \in U_{P_J}, u' \in U_{P^-_{\d(J)}}\}$ of $\Lie(G
\times G)$. We will simply write $h_{J, w_0 w^{\d(J)}_0, \d}$ as
$h_{J, \d}$, $[J, w, v]_{w_0 w^{\d(J)}_0, \d}$ as $[J, w, v]_{\d}$
and $Z^w_{J, w_0 w^{\d(J)}_0, \d}$ as $Z^w_{J, \d}$. We call
$(Z^w_{J, \d})_{J \subset I, w \in W^{\d(J)}}$ the $G$-stable
pieces of $\bar{G}^1$. If $G^1=G$, then $h_{J, id}=h_J$ and $[J,
w, v]_{id}=[J, w, v]$.

Note that $h_J$ corresponds to the $\dim(G)$-dimensional
subalgebra $\{(l u, l u') \mid l \in L_J, u \in U_{P_J}, u' \in
U_{P^-_J}\}$ of $\Lie(G \times G)$. Thus $h_{J, \d}=(1, g_0 \i)
h_{\d(J)}$. Hence $$\eqalignno{& [J, w, v]_{\d}=(B \times B) (\dot
w, \dot v) \cdot h_{J, \d}=(B \times B) (\dot w, \dot v) (1, g_0
\i) \cdot h_J \cr &=(1, g_0 \i) (B \times B) (\dot w, \dot
{\d(v)}) \cdot h_{\d(J)}=(1, g_0 \i) [\d(J), w, \d(v)].}$$

Thus we have the following result.

\proclaim{Proposition} Let $J \subset I$ and $w \in W^{\d(J)}$.
Then the closure of $[J, w, 1]_{\d}$ in $\bar{G^1}$ is $\sqcup_{K
\subset J} \sqcup_{x \in W^{\d(K)}, u \in W_J, \hbox{ and } x \ge
w \d(u)} [K, x, u]_{\d}$.
\endproclaim

\subhead 2.6\endsubhead Define $$\underline{Z_{J, y, \d}}=\{(P,
P', \g) \mid P \in \cp_J, P' \in \cp_{J'}, \g \in \overline{H_{P'}
\backslash A_y(P, P') /H_P}\}$$

with $G \times G$ action defined by $(g_1, g_2) (P, P',
\g)=(^{g_2} P, {}^{g_1} P', g_1 \g g_2 \i)$.

Set $P=P_J$ and $P'={}^{\dot y \i} P_{J'}$. Then $\overline{A_y(P,
P')}$ can be identified with $\overline{L_{P, P', g_0}} g_0$ as
varieties with $P' \times P$ action. Moreover, we have a canonical
isomorphism between $\overline{L_{P, P', g_0}}$ and
$\overline{L_{\d(J)}}$. For $K \subset J$, I will identify
$h_{\d(K)} g_0$ with the corresponding element in
$\overline{A_y(P, P')}$.

Then the $G \times G$-orbits in $\underline{Z_{J, y, \d}}$ are in
one-to-one correspondence with the subsets of $J$, i. e.,
$$\underline{Z_{J, y, \d}}=\sqcup_{K \subset J} (G \times G)
\cdot (P, P', h_{\d(K)} g_0).$$

Set $y_K=y w^{\d(J)}_0 w^{\d(K)}_0$. Note that $U_{P_J} (L_J \cap
U_{P_K})=U_{P_K}$ and $$\eqalignno{& U_{P'} (^{\dot y \i} L_{J'}
\cap U_{P^-_{\d(K)}})=(^{\dot y_K \i} U_{^{\dot y_K} P'}) ^{\dot
y_K \i} (^{\dot y_K \dot y \i} L_{J'} \cap {}^{\dot y_K}
U_{P^-_{\d(K)}}) \cr &={}^{\dot y_K \i} (U_{P_{J'}} (L_{J'} \cap
U_{P_{\Ad(y_K) \d(K)}}))={}^{\dot y_K \i} U_{P_{\Ad(y_K)
\d(K)}}.}$$

The isotropic subgroup of $G \times G$ at $(P, P', h_{\d(K)} g_0)$
is $\{(l_1 u_1, g_0 \i l_2 g_0 u_2) \mid l_1, l_2 \in L_{\d(K)},
l_1 l_2 \i \in Z(L_{\d(K)}), u_1 \in U_{^{\dot y_K \i }
P_{\Ad(y_K) \d(K)}}, u_2 \in U_{P_K}\}$. Now set $Q=P_K,
Q'={}^{\dot y_K \i } P_{\Ad(y_K) \d(K)}$ and $\g=H_{Q'} g_0 H_Q$.
Then $\po(Q', {}^{g_0} Q)=y_K$ and $(Q, Q', \g) \in Z_{K, y_K,
\d}$. The isotropic subgroup of $G \times G$ at $(P, P', h_{\d(K)}
g_0)$ is the same as the isotropic subgroup of $G \times G$ at
$(Q, Q', \g) \in Z_{K, y_K, \d}$. Thus we can identify $(P, P',
h_{\d(K)} g_0)$ with $(Q, Q', \g)$ and $(G \times G) \cdot (P, P',
h_{\d(K)} g_0)$ with $Z_{K, y_K, \d}$ as varieties with $G \times
G$ action. In other words,
$$\underline{Z_{J, y, \d}}=\sqcup_{K \subset J} Z_{K, y w^{\d(J)}_0 w^{\d(K)}_0, \d}.$$

\head 3. Partial order on $\ci_{\d}$
\endhead

In this section, we will only consider subvarieties of $G$ and for
any subvariety $X$ of $G$, we denote by $\bar{X}$ the closure of
$X$ in $G$.

\subhead 3.1 \endsubhead Let $y, w \in W$. Then $y \le w$ if and
only if for any reduced expression $w=s_1 s_2 \cdots s_q$, there
exists a subsequence $i_1<i_2<\cdots<i_r$ of $1, 2, \ldots, q$
such that $y=s_{i_1} s_{i_2} \cdots s_{i_r}$. (see \cite{L2, 2.4})

The following assertion follows from the above property.

(1) If $l(w u)=l(w)+l(u)$, then for any $w_1 \le w$ and $u_1 \le
u$, $w_1 u_1 \le w u$.

(2) Let $u, v \in W$ and $i \in I$. Assume that $s_i v<v$, then $u
\le v \Leftrightarrow s u \le v$.

(3) Let $u, v \in W$ and $i \in I$. Assume that $u<s_i u$, then $u
\le v \Leftrightarrow u \le s_i v$.

The assertion (1) follows directly from the above property. The
proofs of assertions (2) and (3) can be found in \cite{L2, 2.5}.

\subhead 3.2 \endsubhead It is known that $G=\sqcup_{w \in W} B
\dot w B$ and for $w, w' \in W$, $B \dot w B \subset \overline{B
\dot w' B}$ if and only if $w \le w'$. Moreover,
$$\overline{B \dot{s_i} B \dot w B}=\cases \overline{B \dot w B},
& \hbox{ if } s_i w<w; \cr \overline{B \dot s_i \dot w B}, &
\hbox{ if } s_i w>w.
\endcases$$

Similarly, $G=\sqcup_{w \in W} B \dot w B^-$ and for $w, w' \in
W$, $B \dot w B^- \subset \overline{B \dot w' B^-}$ if and only if
$w \ge w'$. Moreover, $$\overline{B \dot{s_i} B \dot w B^-}=\cases
\overline{B \dot w B^-}, & \hbox{ if } s_i w>w; \cr \overline{B
\dot s_i \dot w B^-}, & \hbox{ if } s_i w<w.
\endcases$$

\proclaim{Lemma 3.3} Let $u, w \in W$. Then

(1) The subset $\{v w \mid v \le u\}$ of $W$ contains a unique
minimal element $y$. Moreover, $l(y)=l(w)-l(y w \i)$ and
$\overline{B \dot u B \dot w B^-}=\overline{B \dot y B^-}$.

(2) The subset $\{v w \mid v \le u\}$ of $W$ contains a unique
maximal element $y'$. Moreover, $l(y')=l(w)+l(y' w \i)$ and
$\overline{B \dot u B \dot w B}=\overline{B \dot y' B}$.
\endproclaim

Proof. We will only prove part (1). Part (2) can be proved in the
same way.

For any $v \le u$, $B \dot v \subset \overline{B \dot u B}$. Thus
$\overline{B \dot v \dot w B^-} \subset \overline{B \dot u B} \dot
w B^-\subset \overline{B \dot u B \dot w B^-}$. On the other hand,
$\overline{B \dot u B \dot w B^-}$ is an irreducible, closed, $B
\times B^-$-stable subvariety of $G$. Thus there exists $y \in W$,
such that $\overline{B \dot u B \dot w B^-}=\overline{B \dot y
B^-}$. Since $B \dot v \dot w B^- \subset \overline{B \dot y
B^-}$, we have that $v w \ge y$. Now it suffices to prove that
$y=v w$ for some $v \le u$ with $l(v w)=l(w)-l(v)$.

We argue by induction on $l(u)$. If $l(u)=0$, then $u=1$ and
statement is clear. Assume now that $l(u)>0$. Then there exists $i
\in I$, such that $s_i u<u$. We denote $s_i u$ by $u'$. Now
$$\overline{B \dot u B \dot w B^-}=\overline{B \dot s_i B \dot u'
B \dot w B^-}=\overline{B \dot s_i \overline{B \dot u' B \dot w
B^-}}.$$

By induction hypothesis, there exists $v' \le u'$, such that $l(v'
w)=l(w)-l(v')$ and $\overline{B \dot u' B \dot w B^-}=\overline{B
\dot v' \dot w B^-}$. Thus $$\overline{B \dot s_i \overline{B \dot
u' B \dot w B^-}}=\overline{B \dot s_i \overline{B \dot v' \dot w
B^-}}=\overline{B \dot s_i B \dot v' \dot w B^-}=\cases
\overline{B \dot v' \dot w B^-}, & \hbox{ if } s_i v' w>v' w; \cr
\overline{B \dot s_i \dot v' \dot w B^-}, & \hbox{ if } s_i v'
w<v' w. \endcases$$

Note that $s_i u<u$ and $v' \le s_i u<u$. Thus $s_i v' \le u$.
Moreover, if $s_i v' w<v' w$, then $l(s_i v' w)=l(v'
w)-1=l(w)-l(v')-1$. Thus we have that $l(s_i v')=l(v)+1$ and
$l(s_i v' w)=l(w)-l(s_i v')$. Therefore, the statement holds for
$u$. \qed

\proclaim{Corollary 3.4} Let $u, w, w' \in W$ with $w' \le w$.
Then

(1) There exists $v \le u$, such that $v w' \le u w$.

(2) There exists $v' \le u$, such that $u w' \le v' w$.
\endproclaim

Proof. Let $v \le u$ be the element of $W$ such that $v w'$ is the
unique minimal element in $\{v' w' \mid v' \le u\}$. Then
$\overline{B \dot u B \dot w' B^-}=\overline{B \dot v \dot w'
B^-}$. Since $w' \le w$, we have that $B \dot w B^- \subset
\overline{B \dot w' B^-}$. Thus $$B \dot u \dot w B^- \subset B
\dot u B \dot w B^- \subset B \dot u \overline{B \dot w' B^-}
\subset \overline{B \dot u B \dot w' B^-}=\overline{B \dot v \dot
w' B^-}.$$

So $u w \ge v w'$. Thus Part (1) is proved. Part (2) can be proved
in the same way. \qed

\subhead 3.5\endsubhead We will recall some known results about
$W^J$.

(1) If $w \in W^J$ and $i \in I$, then there are three
possibilities.

$\quad$ (a) $s_i w>w$ and $s_i w \in W^J$;

$\quad$ (b) $s_i w>w$ and $s_i w=w s_j$ for some $j \in J$;

$\quad$ (c) $s_i w<w$ in which case $s_i w \in W^J$.

(2) If $w \in W^J$, $v \in W_J$ and $K \subset J$, then $v \in
W^K$ if and only if $w v \in W^K$.

(3) If $w \in {}^{J'} W^J$ and $u \in W_{J'}$, then $u w \in W^J$
if and only if $u \in W^K$, where $K=J' \cap \Ad(w) J$.

\proclaim{Lemma 3.6} Let $w \in {}^{J'} W^J$, $u \in W_{J'}$ and
$K=J' \cap \Ad(w) J$, then $u w=v w u'$ for some $v \in W_{J'}
\cap W^K$ and $u' \in W_{\Ad(w \i) K}$.
\endproclaim

Proof. We argue by induction on $l(u)$. If $u=1$, then the
statement is clear. Now assume that $u=s_i u_1$ for some $i \in
J'$ and $l(u_1)<l(u)$. Then by induction hypothesis, $u_1 w=v_1 w
u'_1$ for some $v_1 \in W_{J'} \cap W^K$ and $u'_1 \in W_{\Ad(w
\i) K}$.

If $s_i v_1 w \in W^J$, then the statement holds for $u$. Now
assume that $s_i v_1 w \notin W^J$. Then $s_i v_1 w>v_1 w$. Hence
$s_i v_1>v_1$. Moreover, $s_i v_1 \notin W^K$. Thus $s_i v_1=v_1
s_k$ for some $k \in K$. Note that $s_k w=w s_l$ for some $l \in
\Ad(w \i) K$. Thus the statement holds for $u$. The lemma is
proved. \qed

\subhead 3.7\endsubhead Let $J \subset I$ and $w, w' \in W$ with
$l(w)=l(w')$. We say that $w'$ can be obtained from $w$ via a $(J,
\d)$-cyclic shift if $w=s_{i_1} s_{i_2} \cdots s_{i_n}$ is a
reduced expression and either (1) $i_1 \in J$ and $w'=s_{i_1} w
s_{\d(i_1)}$ or (2) $i_n \in \d(J)$ and $w'=s_{\d \i(i_n)} w
s_{i_n}$. We say that $w$ and $w'$ are equivalent in J if there
exists a finite sequences of elements $w=w_0, w_1, \ldots, w_m=w'$
such that $w_{k+1}$ can be obtained from $w_k$ via a $(J,
\d)$-cyclic shift. (We then write $w \sim_{J, \d} w'$.)

\proclaim{Proposition 3.8} Let $(J, w) \in \ci_{\d}$ and $w' \in
W$. The following conditions on $w'$ are equivalent:

(1) $w' \ge u \i w \d(u)$ for some $u \in W_J$.

(2) $w' \ge u \i w \d(v)$ for some $u \le v \in W_J$.

(3) $w' \ge w_1$ for some $w_1 \sim_{J, \d} w$.
\endproclaim

Proof. The implication (1)$\Rightarrow$(2) is trivial. The
implication (3)$\Rightarrow$(1) follows from the definition. We
now prove the implication (2)$\Rightarrow$(3) by induction on
$|J|$. Assume that the implication holds for all $J' \subset I$
with $|J'|<|J|$. Then we prove that the implication holds for $J$
by induction on $l(v)$.

Set $w=x y$ with $x \in W_J$ and $y \in {}^J W^{\d(J)}$. Set $K=J
\cap \d\i\Ad(y \i) J$, $v=v_1 v_2$ with $v_1 \in W_K, v_2 \in {}^K
W$ and $u=u_1 u_2$ with $u_1 \le v_1$, $u_2 \le v_2$ and
$l(u)=l(u_1)+l(u_2)$. There are two cases.

Case 1. $u_2=v_2=1$.

In this case, $u, v \in W_K$ and $w \in W^{\d(K)}$. If $|K|<|J|$,
then by induction hypothesis, $u \i w \d(v) \ge w_1$ for some $w_1
\sim_{K, \d} w$. If $K=J$, then since $w=x y \in W^{\d(J)}$, we
have that $x=1$. Thus $u \i w \d(v) \ge w$. The implication is
proved in this case.

Case 2. $v_2 \neq 1$.

In this case, $l(v_1)<l(v)$. By induction hypothesis, there exists
$w_1 \sim_{J, \d} w$, such that $w_1 \le u_1 \i w \d(v_1)$. Let
$u_3 \le u_2$ be the element in $W$ such that $u_3 \i w_1$ is the
unique minimal element in $\{(u') \i w_1 \mid u' \le u_2\}$. Then
$l(u_3 \i w_1)=l(w_1)-l(u_3)$ and $u_3 \i w_1 \le u_2 \i u_1 \i w
\d(v_1)=u \i w \d(v_1)$. By 3.6, $u \i w=a b$ for some $a \in
W^{\d(J)}$ and $b \in W_{\d(K)}$. Thus $l(u \i w \d(v_1 v_2))=l(a
b \d(v_1 v_2))=l(a)+l(b \d(v_1 v_2))=l(a)+l(b
\d(v_1))+l(\d(v_2))=l(a b \d(v_1))+l(\d(v_2))=l(u \i w
\d(v_1))+l(\d(v_2))$. By 3.1, $u_3 \i w_1 \d(u_3) \le u \i w
\d(v)$.

Now assume that $u_3=s_{i_1} s_{i_2} \cdots s_{i_k}$ and $u_3 \i
w_1=s_{j_1} s_{j_2} \cdots s_{j_l}$ are reduced expressions. For
$m=1, 2, \ldots, k+1$, set $$x_m=(s_{i_m} s_{i_{m+1}} \cdots
s_{i_k}) (s_{j_1} s_{j_2} \cdots s_{j_l}) (s_{\d(i_1)} s_{\d(i_2)}
\cdots s_{\d(i_{m-1})}).$$ Then $l(x_m) \le k+l=l(w_1)$ for all
$m$. On the other hand, for any $m$, there exists $y_m \in W_J$,
such that $x_m=y_m \i w \d(y_m)$. Note that $w \in W^J$, we have
that $l(y_m \i w \d(y_m)) \ge l(w \d(y_m))-l(y_m \i)=l(w)=l(w_1)$
for all $y_m \in W_J$. Thus $l(x_m)=l(w_1)$ and $x_m \sim_{J, \d}
w_1$ for all $m$. In particular, $u_3 \i w_1 \d(u_3)=x_{k+1}
\sim_{J, \d} w_1$. The implication is proved in this case. \qed

\subhead Remark \endsubhead We see from the proof that $u \i w
\d(v) \ge x \i w \d(x)$ for some $x \le u$. This result will be
used in the proof of 5.2.

\subhead 3.9\endsubhead Let $(J, w) \in \ci_{\d}$ and $w' \in W$,
we say that $w' \ge_{J, \d} w$ if $w'$ satisfies the equivalent
conditions 3.8 (1)-(3). It is easy to see that $x \ge w
\Rightarrow x \ge_{J, \d} w \Rightarrow l(x) \ge l(w)$.

Now for $(J_1, w_1), (J_2, w_2) \in \ci_{\d}$, we say that $(J_1,
w_1) \le_{\d} (J_2, w_2)$ if $J_1 \subset J_2$ and $w_1 \ge_{J,
\d} w_2$. In the end of this section, we will show that $\le$ is a
partial order on $\ci_{\d}$. (The definition of partial order can
be found in 3.12).

\proclaim{Lemma 3.10} Let $J \subset I$, $w \in W^J$, $u \in W$
with $l(u w)=l(u)+l(w)$. Assume that $u w=x v$ with $x \in W^J$
and $v \in W_J$. Then for any $v' \le v$, there exists $u' \le u$,
such that $u' w=x v'$.
\endproclaim

Proof. We argue by induction on $l(u)$. If $l(u)=0$, then $u=1$
and statement is clear. Assume now that $l(u)>0$. Then there
exists $i \in I$, such that $s_i u<u$. We denote $s_i u$ by $u_1$.
Let $u_1 w=x_1 v_1$ with $x_1 \in W^J$ and $v_1 \in W_J$. Then
$s_i x_1>x_1$.

If $s_i x_1 \in W^J$, then the lemma holds by induction
hypothesis. If $s_i x_1 \notin W^J$, then there exists $j \in J$,
such that $s_i x_1=x_1 s_j$. In this case, $s_j v_1>v_1$. Let $v'
\le s_j v_1$. If $v' \le v_1$, then the lemma holds by induction
hypothesis. If $v' \nleqslant v_1$, then $v'=s_j v'_1$ for some
$v'_1 \le v_1$. By induction hypothesis, there exists $u'_1 \le
u_1$, such that $u'_1 w=x_1 v'_1$. Thus $s_i u'_1 w=x_1 s_j v'_1$.
The lemma holds in this case. \qed

\proclaim{Lemma 3.11} Fix $J \subset I$ and $w \in W^{\d(J)}$. For
any $K \subset J$, $w' \in W^{\d(K)}$ with $w' \ge_{J, \d} w$,
there exists $x \in W^{\d(K)}$, $u \in W_J$ and $u_1 \in W_K$,
such that $x \ge w \d(u)$ and $w'=u_1 \i u \i x \d(u_1)$.
\endproclaim

Proof. Since $w' \ge_{J, \d} w$, there exists $v_1 \in W_J$, such
that $w' \ge v_1 \i w \d(v_1)$. By 3.4, there exists $v' \le v_1$,
such that $v' w' \ge w \d(v_1) \ge w \d(v')$. Let $v$ be a minimal
element in the set $\{v \in W_J \mid v w' \ge w \d(v)\}$. Then
$l(v w')=l(v)+l(w')$. Now assume that $v w'=x \d(v')$ for some $x
\in W^{\d(K)}$ and $v' \in W_K$. Then there exists $v'_1 \le v'$,
such that $x \ge w \d(v) \d(v'_1) \i$. By 3.10, $x \d(v'_1)=v_2
w'$ for some $v_2 \le v$. Since $l(x \d(v'_1))=l(x)+l(v'_1)$, $v_2
w'=x \d(v'_1) \ge w \d(v) \ge w \d(v_2)$. Therefore, $v_2=v$ and
$v'_1=v'$. So $x \ge w \d(v) \d(v') \i$. Now set $u=v (v') \i$ and
$u_1=v'$. Then $w'=v \i x \d(v')=u_1 \i u \i x \d(u_1)$. \qed

\subhead 3.12\endsubhead A relation $\le$ is a partial order on a
set $S$ if it has:

1. Reflexivity: $a \le a$ for all $a \in S$.

2. Antisymmetry: $a \le b$ and $b \le a$ implies $a=b$.

3. Transitivity: $a \le b$ and $b \le c$ implies $a \le c$.

\proclaim{Proposition 3.13} The relation $\le_{\d}$ on the set
$\ci_{\d}$ is a partial order.
\endproclaim

Proof. Reflexivity is clear from the definition.

For $(J_1, w_1), (J_2, w_2) \in \ci_{\d}$ with $(J_1, w_1)
\le_{\d} (J_2, w_2)$ and $(J_2, w_2) \le_{\d} (J_1, w_1)$, we have
that $J_1=J_2$ and $l(w_1)=l(w_2)$. Since $w_1 \ge w'_2$ for some
$w'_2 \sim_{J, \d} w_2$ and $l(w_1)=l(w_2)=l(w'_2)$, $w_1=w'_2 \in
W^{\d(J_2)}$. Hence $w_1=w'_2=w_2$. Therefore $(J_1, w_1)=(J_2,
w_2)$. Antisymmetry is proved.

Let $(J_1, w_1), (J_2, w_2)$ and $(J_3, w_3) \in \ci_{\d}$. Assume
that $(J_1, w_1) \le_{\d} (J_2, w_2)$ and $(J_2, w_2) \le_{\d}
(J_3, w_3)$. Then $J_1 \subset J_2 \subset J_3$. Moreover, there
exists $x \in W^{\d(J_2)}$, $u \in W_{J_3}$ and $u_1 \in W_{J_2}$,
such that $x \ge w_3 \d(u)$ and $w_2=u_1 \i u \i x \d(u_1)$. Since
$w_1 \ge_{J_2, \d} w_2$, there exists $u_2 \in W_{J_2}$, such that
$w_1 \ge u_2 \i u \i x \d(u_2)$. Note that $l(x
\d(u_2))=l(x)+l(u_2)$ and $x \ge w_3 \d(u)$. Thus $x \d(u_2) \ge
w_3 \d(u u_2)$. By 3.4, there exists $v \le u u_2$, such that $w_1
\ge v \i w_3 \d(u u_2)$. By 3.7, $w_1 \ge_{J_3, \d} w_3$.
Transitivity is proved. \qed

\head 4. The closure of any $G$-stable piece \endhead

\subhead 4.1\endsubhead We have that $$G^1=\sqcup_{w \in W} B \dot
w U^- \dot w^{\d(J)}_0 g_0=\sqcup_{w \in W} B \dot w \dot
w^{\d(J)}_0 U_{P^-_{\d(J)}} U_{\d(J)} g_0.$$

Moreover, $B \dot w U^-=\sqcup_{b \in U^-_J \cap {}^{\dot w \i}
U^-} B \dot w U_{P^-_J} b$. Thus $$\eqalignno{B \dot w U_{P^-_J}
U_J &=B \dot w \dot w^J_0 U^- \dot w^J_0=\sqcup_{b \in U^-_J \cap
{}^{(\dot w \dot w^J_0) \i} U^-} B \dot w \dot w^J_0 U_{P^-_J} b
\dot w^J_0 \cr &=\sqcup_{b \in U_J \cap {}^{\dot w \i} U^-} B \dot
w U_{P^-_J} b.}$$

Note that if $w=w' u$ with $w' \in W^J$ and $u \in W_J$, then
$$U_J \cap {}^{\dot w \i} U^-={}^{\dot u \i} (^{\dot u} U_J \cap {}^{(\dot w') \i} U^-)={}^{\dot u \i} (^{\dot u} U_J \cap U^-_J)=U_J
\cap {}^{\dot u \i} U^-_J.$$

\proclaim{Lemma 4.2} Let $(J, w) \in \ci_{\d}$. For any $u \in W$
and $b \in B$, there exists $v \le u$, such that $\dot u b \dot w
\in B \dot v \dot w U_{P^-_{\d(J)}} U_{\d(J)}$.
\endproclaim

Proof. We will prove the statement by induction on $l(u)$.

If $u=1$, then the statement holds. If $u=s_i u_1$ with
$l(u_1)=l(u)-1$, then by induction hypothesis, there exists $v_1
\le u_1$, such that $\dot u_1 b \dot w \in b' \dot v_1 \dot w
U_{P^-_{\d(J)}} U_{\d(J)}$ for some $b' \in B$. Write $b'=b_1
b_2$, where $b_1 \in U_{P_{\{i\}}}$ and $b_2 \in U_{\{i\}}$. Then
$\dot s_i b' \dot v_1 \dot w=(\dot s_i b_1 \dot s_i \i) \dot s_i
b_2 \dot v_1 \dot w$ with $\dot s_i b_1 \dot s_i \i \in B$.

If $(\dot v_1 \dot w) \i b_2 \dot v_1 \dot w \in U_{P^-_{\d(J)}}
U_{\d(J)}$, then $\dot s_i b_2 \dot v_1 \dot w \in \dot s_i \dot
v_1 \dot w U_{P^-_{\d(J)}} U_{\d(J)}$. Otherwise, $b_2 \neq 1$ and
$(\dot v_1 \dot w) \i U^-_{\{i\}} \dot v_1 \dot w \subset
U_{P^-_{\d(J)}} U_{\d(J)}$. Note that $\dot s_i b_2 \in B
U^-_{\{i\}}$. Thus $\dot s_i b_2 \dot v_1 \dot w \in B \dot v_1
\dot w U_{P^-_{\d(J)}} U_{\d(J)}$. The statement holds in both
cases. \qed

\subhead 4.3 \endsubhead Let $z \in (G, 1) \cdot h_{J, \d}$. Then
$z$ can be written as $z=(b \dot w \dot u b', 1) \cdot h_{J, \d}$
with $b \in B$, $w \in W^{\d(J)}$, $u \in W_{\d(J)}$ and $b' \in
U_{\d(J)} \cap {}^{\dot u \i} U^-_{\d(J)}$. Moreover, $w, u, b'$
are uniquely determined by $z$.

Set $J_0=J$. To $z \in (G, 1) \cdot h_{J, \d}$, we associate a
sequence $(J_i, w_i, v_i, v'_i, c_i, z_i)_{i \ge 1}$ with $J_i
\subset J$, $w_i \in W^{\d(J)}$, $v_i \in W_{J_{i-1}} \cap
{}^{J_i} W$, $v'_i \in W_{J_i}$, $c_i \in U_{\d(J_{i-1})} \cap
{}^{\dot \d(v_i \i)} U^-_{\d(J_{i-1})}$ and $z_i \in (B \dot w_i
\dot \d(v'_i) U_{\d(J)} \dot \d(v_i) c_i, 1) \cdot h_{J, \d}$ and
in the same $G_{diag}$-orbit as $z$. The sequence is defined as
follows.

Assume that $z \in (B \dot w \dot \d(u) U_{\d(J)}, 1) \cdot h_{J,
\d}$ with $w \in W^{\d(J)}$ and $u \in W_{\d(J)}$. Then set
$z_1=z$, $J_1=J$, $w_1=w$, $v_1=1$, $v'_1=u$ and $c_1=1$.

Assume that $k \ge 1$, that $J_k, w_k, v_k, v'_k, c_k, z_k$ are
already defined and that $J_k \subset J_{k-1}$, $w_k \in
W^{\d(J)}$, $W_{J_{k-1}} w_k \subset W^{\d(J)} W_{\d(J_k)}$, $v_k
\in W_{J_{k-1}} \cap {}^{J_k} W$, $v'_k \in W_{J_k}$, $c_k \in
U_{\d(J_{k-1})} \cap {}^{\dot \d(v_k \i)} U^-_{\d(J_{k-1})}$ and
$z_k \in (B \dot w_k \dot \d(v'_k) U_{\d(J)} \dot \d(v_k) c_k, 1)
\cdot h_{J, \d}$. Set $z_{k+1}=(g_0 \i \dot \d(v_k) c_k g_0, g_0
\i \dot \d(v_k) c_k g_0) z_k$. Then $z_{k+1} \in (G, 1) \cdot
h_{J, \d}$. Moreover, by 4.2, there exists $x_k \le v_k$, such
that $z_{k+1} \in (B \dot x_k w_k \dot \d(v'_k) U_{\d(J)}, 1)
\cdot h_{J, \d}$.

Let $y_{k+1}$ be the unique element of the minimal length in
$W_{J_k} x_k w_k \d(v'_k) W_{\d(J)}$. Set $J_{k+1}=J_k \cap \d \i
\Ad(y_{k+1} \i) J_k$. Since $W_{J_{k-1}} w_k \subset W^{\d(J)}
W_{\d(J_k)}$, then we have that $x_k w_k \d(v'_k)=w_{k+1}
\d(v'_{k+1} v_{k+1})$ for some $w_{k+1} \in W^{\d(J)}$, $v'_{k+1}
\in W_{J_{k+1}}$ and $v_{k+1} \in W_{J_k} \cap {}^{J_{k+1}} W$.
Note that $W_{J_k} w_{k+1} \subset W^{\d(J)} W_{\d(J) \cap
\Ad(y_{k+1} \i) J_k}$. On the other hand, $W_{J_k} w_{k+1} \subset
W_{J_{k-1}} w_k W_{\d(J_k)} \subset W^{\d(J)} W_{\d(J_k)}$. Thus
$W_{J_k} w_{k+1} \subset (W^{\d(J)} W_{\d(J) \cap \Ad(y_{k+1} \i)
J_k}) \cap (W^{\d(J)} W_{\d(J_k)})=W^{\d(J)} W_{\d(J_{k+1})}$.

Moreover $z_{k+1} \in (B \dot w_{k+1} \dot \d(v'_{k+1}) U_{\d(J)}
\dot \d(v_{k+1}) c_{k+1}, 1) \cdot h_{J, \d}$ for a unique
$c_{k+1} \in U_{J_k} \cap {}^{\dot \d(v_{k+1} \i)} U^-_{J_k}$.

This completes the inductive definition. Moreover, for sufficient
large $n$, we have that $J_n=J_{n+1}=\cdots$,
$w_n=w_{n+1}=\cdots$, $v'_n=v'_{n+1}=\cdots$ and
$v_n=v_{n+1}=\cdots=1$.

\subhead 4.4\endsubhead Let $K \subset J$, $y \in {}^K W^{\d(K)}$
and $K=\Ad(y) \d(K)$. Then for any $u \in W_K$, we have that $$(B
\dot y \dot \d(u) U_{\d(K)}, 1) \cdot h_{J, \d} \subset G_{diag}
(\dot y L_{\d(K)}, B) \cdot h_{J, \d}=G_{diag} (\dot y L_{\d(K)},
U_{P_K}) \cdot h_{J, \d}.$$

For any $l \in L_K$, there exists $l' \in L_K$, such that $l' \dot
y g_0 l (l') \i \in \dot y g_0 (L_K \cap B)$. Thus $(L_K)_{diag}
(\dot y (L_{\d(K)} \cap B), U_{P_K}) \cdot h_{J, \d}=(\dot y
L_{\d(K)}, U_{P_K}) \cdot h_{J, \d}$. Hence $$(B \dot y \dot \d(u)
U_{\d(K)}, 1) \cdot h_{J, \d} \subset G_{diag} (\dot y (L_{\d(K)}
\cap B), U_{P_K}) \cdot h_{J, \d}=G_{diag} (\dot y, B) \cdot h_{J,
\d}=Z^y_{J, \d}.$$

Now for any $z \in (G, 1) \cdot h_{J, \d}$, let $(z_i, J_i, w_i,
v_i, v'_i, c_i)_{i \ge 1}$ be the sequence associated to $z$.
Assume that $J_n=J_{n+1}=\cdots$, $w_n=w_{n+1}=\cdots$,
$v'_n=v'_{n+1}=\cdots$ and $v_n=v_{n+1}=\cdots=1$. Then we have
showed that $z_n \in Z^{w_n}_{J, \d}$. Thus $z \in Z^{w_n}_{J,
\d}$.

Note that for any $z \in Z_{J, \d}$, $z$ is in the same $G$-orbit
as an element of the form $(G, 1) \cdot h_{J, \d}$. Therefore,
given $z \in Z_{J, \d}$, our procedure determines the $G$-stable
piece $Z^w_{J, \d}$ that contains $z$.

${ }$

Now we are able to describe the closure of $Z^w_{J, \d}$. In 4.5,
we will only consider subvarieties of $\bar{G^1}$ and for any
subvariety $X$ of $\bar{G^1}$, we denote by $\bar{X}$ the closure
of $X$ in $\bar{G^1}$.

\proclaim{Theorem 4.5} For any $(J, w) \in \ci_{\d}$, we have that
$$\overline{Z^w_{J, \d}}=\sqcup_{(K, w') \le_{\d} (J, w)} Z^{w'}_{K, \d}.$$
\endproclaim

Proof. Define $\pi': G \times \overline{[I, 1,
1]_{\d}}@>>>\bar{G^1}$ by $\pi(g, z)=(g, g) \cdot z$. The morphism
is invariant under the $B$-action defined by $b (g, z)=(g b \i,
\pi'(b, z))$. Denote by $G \times_B \overline{[I, 1, 1]_{\d}}$ the
quotient, we obtain a morphism $\pi: G \times_B \overline{[I, 1,
1]_{\d}}@>>>\bar{G^1}$. Because $G/B$ is projective, $\pi$ is
proper and hence surjective.

Note that $\overline{[J, w, 1]_{\d}}=\sqcup_{K \subset J}
\sqcup_{x \in W^{\d(K)}, u \in W_J, \hbox{ and } x \ge w \d(u)}
[K, x, u]_{\d}$. Since $Z^w_{J, \d}=\pi(G \times_B [J, w,
1]_{\d})$, we have that
$$\overline{Z^w_{J, \d}}=\sqcup_{K \subset J} \cup_{x \in W^{\d(K)}, u \in
W_J, \hbox{ and } x \ge w \d(u)} G_{diag} \cdot [K, x, u]_{\d}.$$

For any $z \in [K, x, u]_{\d}$ with $x \in W^{\d(K)}$, $u \in W_J$
and $x \ge w \d(u)$, we have that $z \in (B \dot x, B \dot u)
\cdot h_{K, \d}=G_{diag} (\dot u \i B \dot x, 1) \cdot h_{K, \d}
\subset \sqcup_{v \le u \i} G_{diag} (B \dot v \dot x U_{\d(K)},
1) \cdot h_{K, \d}$.

Fix $v \le u \i$ and $z' \in (B \dot v \dot x U_{\d(K)}, 1) \cdot
h_{K, \d}$. Let $(z_i, J_i, w_i, v_i, v'_i, c_i)_{i \ge 1}$ be the
sequence associated to $z'$. Then for any $i$, there exists $x_i
\le v_i$, such that $x_i w_i \d(v'_i)=w_{i+1} \d(v'_{i+1}
v_{i+1})$. Assume that $J_n=J_{n+1}=\cdots$, $w_n=w_{n+1}=\cdots$,
$v'_n=v'_{n+1}=\cdots$ and $v_n=v_{n+1}=\cdots=1$. Set
$x_{\infty}=x_n x_{n-1} \cdots x_2$ and $v_{\infty}=v'_n (v_n
v_{n-1} \cdots v_2)$. Note that $x_1=v_1=1$. Then $x_{\infty} v
x=x_{\infty} w_1 \d(v'_1)=w_n \d(v_{\infty})$. Since $v'_n \in
W_{J_{n+1}}$ and $v_i \in W_{J_i} \cap {}^{J_{i+1}} W$, we have
that $l(v_{\infty})=l(v'_n)+l(v_n)+l(v_{n-1})+\cdots+l(v_2)$. Thus
$x_{\infty} \le v_{\infty}$. By 4.4, $z' \in Z_{K, x_{\infty} v x
\d(v_{\infty} \i)}$. Note that $v \i \le u$ and $l(w
u)=l(w)+l(u)$. Thus $w \d(v \i) \le w \d(u) \le x$. Similarly, $w
\d(v \i x_{\infty} \i) \le x \d(v_{\infty} \i)$. By 3.4, there
exist $v' \le v \i x_{\infty} \i$, such that $(v') \i w \d(v \i
x_{\infty} \i) \le x_{\infty} v x \d(v_{\infty} \i)$. Thus by 3.8,
$x_{\infty} v x \d(v_{\infty} \i) \ge_{J, \d} w$.

For any $K \subset J$ and $w' \in W^{\d(K)}$ with $w' \ge_{J, \d}
w$, there exists $x \in W^{\d(K)}$, $u \in W_J$ and $u_1 \in W_K$,
such that $x \ge w \d(u)$ and $w'=u_1 \i u \i x \d(u_1)$. Since
$[K, x, u]_{\d} \subset \overline{[J, w, 1]_{\d}}$. We have that
$(\dot x T, \dot u) \cdot h_{K, \d} \subset \overline{Z^w_{J,
\d}}$. Therefore $(\dot u \i \dot x T, 1) \cdot h_{K, \d} \subset
\overline{Z^w_{J, \d}}$. Note that $u \i x=u_1 w' \d(u_1) \i$.
Then $(\dot u_1 \dot w' \dot \d(u_1) \i T, 1) \cdot h_{K, \d}
\subset \overline{Z^w_{J, \d}}$. Thus $$(\dot u_1 \i, \dot u_1 \i)
(\dot u_1 \dot w' \dot \d(u_1) \i T, 1) \cdot h_{K, \d}=(\dot w'
T, 1) \cdot h_{K, \d} \subset \overline{Z^w_{J, \d}}.$$

By 1.12, $Z^{w'}_{K, \d} \subset \overline{Z^w_{J, \d}}$. The
theorem is proved. \qed

${}$

Our method also works in another situation.

\proclaim{Proposition 4.6} The closure of $Z^w_{J, 1, \d}$ in
$Z_{J, 1, \d}$ is $\sqcup_{w' \in W^{\d(J)}, w \ge_J w'}
Z^{w'}_{J, 1, \d}$.
\endproclaim

Proof. In the proof, we will only consider subvarieties of $Z_{J,
1, \d}$ and for any subvariety $X$ of $Z_{J, 1, \d}$, we denote by
$\bar{X}$ its closure in $Z_{J, 1, \d}$.

Note that the morphism $\pi: Z_{J, 1, \d}@>>>\cp_J$ defined by
$\pi(P, Q, \g)=P$ for $(P, Q, \g) \in Z_{J, 1, \d}$ is a locally
trivial fibration with isomorphic fibers. Moreover, $i: \pi
\i(P_J)@>>>G^1/H_{P_J}$ defined by $i(P, Q, \g)=\g$ for $(P, Q,
\g) \in \pi \i(P_J)$ is an isomorphism. Now $[J, w, 1]_{1, \d}
\subset \pi \i(P_J)$ and $i([J, w, 1]_{1, \d})=B \dot w B g_0/
H_{P_J}$. Thus $\overline{[J, w, 1]_{1, \d}}=\sqcup_{w' \le w} [J,
w', 1]_{1, \d}$. For any $w' \in W^J$ with $w \ge_{J, \d} w'$,
there exists $u \in W_J$, such that $w \ge u \i w' \d(u)$. Thus
$(\dot u \i \dot w' \dot \d(u) T, 1) \cdot h_{J, 1, \d} \subset
\overline{[J, w, 1]_{1, \d}}$. Hence $G_{diag} (\dot u \i \dot w'
\dot \d(u) T, 1) \cdot h_{J, 1, \d}=G_{diag} (\dot w' T, 1) \cdot
h_{J, 1, \d} \subset \overline{Z^w_{J, 1, \d}}$. So $Z^{w'}_{J, 1,
\d} \subset \overline{Z^w_{J, 1, \d}}$.

On the other hand, for any $z \in [J, w', 1]_{1, \d}$, by the
similar argument as we did in 4.3 and 4.4, there exists $u \le v
\in W_J$, such that $u w' \d(v \i) \in W^{\d(J)}$ and $z \in Z^{u
w' v \i}_{J, 1, \d}$. If moreover, $w' \le w$, then $w \ge u \i (u
w' \d(v \i)) \d(v)$. Thus $w \ge_{J, \d} u w' v \i$. Therefore $z
\in \sqcup_{w' \in W^{\d(J)}, w \ge_{J, \d} w'} Z^{w'}_{J, 1,
\d}$. The proposition is proved. \qed

\head 5. The cellular decomposition \endhead

\subhead 5.1\endsubhead A finite partition of a variety $X$ into
subsets is said to be an $\a$-partition if the subsets in the
partition can be indexed $X_1, X_2, \ldots, X_n$ in such a way
that $X_1 \cup X_2 \cup \cdots \cup X_i$ is closed in $X$ for
$i=1, 2, \ldots, n$. We say that a variety has a cellular
decomposition if it admits an $\a$-partition into subvarieties
which are affine spaces. It is easy to see that if a variety $X$
admits an $\a$-partition into subvarieties and each subvariety has
a cellular decomposition, then $X$ has a cellular decomposition.

\proclaim{Lemma 5.2} Let $(J, w) \in \ci_{\d}$, $K \subset J$ and
$w' \in W$ with $\Ad(w') \d(K)=K$. If $w' v \ge_{J, \d} w$ for
some $v \in W_{\d(K)}$, then $w' \ge_{J, \d} w$.
\endproclaim

Proof. Fix $w'$ and $(J, w)$. It suffices to prove the following
statement:

Let $u \in W_J$ and $v \in W_{\d(K)}$. If $w' v \ge u \i w \d(u)$,
then $w' \ge_{J, \d} w$.

We argue by induction on $l(u)$. Assume that the statement holds
for all $u'<u$. Then I will prove that the statement holds for $u$
by induction on $l(v)$. If $l(v)=0$, then $v=1$ and the statement
holds in this case. Now assume that $l(v)>0$.

Set $u=u_1 u_2$ with $u_1 \in W^K$ and $u_2 \in W_K$. If $u_2=1$,
then $u \in W^K$ and $w \d(u) \in W^{\d(K)}$. By 3.4, there exists
$u' \le u$, such that $u' w' v \ge w \d(u)$. Assume that $v=v'
s_k$ for $v'<v$ and $k \in \d(K)$. Then $w \d(u)<w \d(u) s_k$. By
3.1, $w \d(u) \le u' w' v'$. By 3.4, there exists $u'_1 \le u' \le
u$, such that $w' v' \ge (u'_1) \i w \d(u)$. By the remark of 3.8,
$w' v' \ge (u'_2) \i w \d(u'_2)$ for some $u'_2 \le u'_1$. Thus by
induction hypothesis, $w' \ge_{J, \d} w$.

If $u_2 \neq 1$. Then $l(u_1)<l(u)$. By 3.4, there exists $u_3 \le
u_2$ and $u_4 \le u_2 \i$, such that $u_3 w' v \d(u_4) \ge u_1 \i
w \d(u_1)$. Note that $u_3 w' v u_4=w' ((w') \i u_3 w') v \d(u_4)
\in w' W_{\d(K)}$. By induction hypothesis on $l(u_1)$, $w'
\ge_{J, \d} w$. \qed

\subhead 5.3\endsubhead Let $J \subset I$. For $w \in W$, set
$$\eqalignno{I_1(J, w, \d) &=\max\{K \subset J \mid w \in W^{\d(K)}\}, \cr
I_2(J, w, \d) &=\max\{K \subset J \mid \Ad(w)
\Phi_{\d(K)}=\Phi_K\}.}$$

Now let $(J, w) \in \ci_{\d}$. Set $$W_{\d}(J, w)=\{u \in W \mid u
\ge_{J, \d} w, I_2(J, u, \d) \subset I_1(J, u, \d)\}.$$

For any $u \in W_{\d}(J, w)$, set $$\eqalignno{X^{(J, w, \d)}_u
&=\sqcup_{K \subset I_1(J, u, \d)} \sqcup_{v \in W_{\d(I_2(J, u,
\d))} \cap W^{\d(K)}} Z^{u v}_{K, \d} \cr &=\sqcup_{v \in
W_{\d(I_2(J, u, \d))}} \sqcup_{K \subset I_1(J, u v, \d)} Z^{u
v}_{K, \d}.}$$

For $w' \ge_{J, \d} w$, we have that $w'=u v$ for some $u \in
W^{\d(I_2(J, w'))}$ and $v \in W_{\d(I_2(J, w'))}$. Then $I_2(J,
u, \d)=I_2(J, w', \d) \subset I_1(J, u, \d)$. By 5.2, $u \ge_{J,
\d} w$. Thus $u \in W_{\d}(J, w)$ and $\sqcup_{K \subset I_1(J,
w', \d)} Z^{w'}_{K, \d} \subset X^{(J, w, \d)}_u$.

For $u_1, u_2 \in W(J, w)$ and $v_1 \in W_{\d(I_2(J, u_1))}, v_2
\in W_{\d(I_2(J, u_2))}$ with $u_1 v_1=u_2 v_2$, we have that
$I_2(J, u_1, \d)=I_2(J, u_1 v_1, \d)=I_2(J, u_2 v_2, \d)=I_2(J,
u_2, \d)$. Note that $u_1, u_2 \in W^{\d(I_2(J, u_1))}$. Thus
$u_1=u_2$ and $v_1=v_2$.

Therefore $\overline{Z^w_{J, \d}}=\sqcup_{u \in W_{\d}(J, w)}
X^{(J, w, \d)}_u$.

\proclaim{Lemma 5.4} Let $(J, w) \in \ci_{\d}$. Set $I_2=I_2(J, w,
\d)$. For $K \subset J$, we have that
$$\sqcup_{v \in W_{\d(I_2)} \cap W^{\d(K)}}
(L_{I_2})_{diag} (\dot w \dot v, B \cap L_{I_2}) \cdot h_{K,
\d}=(L_{I_2}, L_{I_2}) (\dot w, 1) \cdot h_{K, \d}.$$
\endproclaim

Proof. At first, we will prove the case when $K \subset I_2$. In
this case, set $g_1=g_0 \dot w$. Then $g_1 L_{\d(I_2)} g_1
\i=L_{\d(I_2)}$ and $g_1 (L_{\d(I_2)} \cap B) g_1 \i=L_{\d(I_2)}
\cap B$. Now consider $\overline{L_{\d(I_2)}/Z(L_{\d(I_2)})} g_1$
(a variety that is isomorphic to
$\overline{L_{\d(I_2)}/Z(L_{\d(I_2)})}$, but with ``twisted''
$L_{\d(I_2)} \times L_{\d(I_2)}$ action, see 2.3). We have that
$$\sqcup_{v \in W_{\d(I_2)} \cap W^{\d(K)}}
(L_{\d(I_2)})_{diag} (\dot v, B \cap L_{\d(I_2)}) \cdot (h_{\d(K)}
g_1)=(L_{\d(I_2)}, L_{\d(I_2)}) \cdot (h_{\d(K)} g_1).$$

(In the case when $g_1^n \in L_{\d(I_2)}$ for some $n \in \bold
N$, $L_{\d(I_2)} g_1$ is a connected component of the group
generated by $L_{\d(I_2)}$ and $g_1$. In this case, the left hand
side is the union of some $L_{\d(I_2)}$-stable pieces and the
equality follows from \cite{L4, 12,3}. The general case can be
shown in the same way.)

Therefore $$\eqalignno{& \sqcup_{v \in W_{\d(I_2)} \cap W^{\d(K)}}
(\dot w \i, g_0)(L_{I_2})_{diag} (\dot w, g_0 \i)(\dot v, B \cap
L_{\d(I_2)}) \cdot h_{\d(K)} \cr &=\sqcup_{v \in W_{\d(I_2)} \cap
W^{\d(K)}} (1, g_1)(L_{\d(I_2)})_{diag} (1, g_1 \i)(\dot v, B \cap
L_{\d(I_2)}) \cdot h_{\d(K)} \cr &=(L_{\d(I_2)}, L_{\d(I_2)})
(\dot w, 1) \cdot h_{\d(K)}.}$$

Note that $h_{K, \d}=h_{\d(K)} g_0$. Then $\sqcup_{v \in
W_{\d(I_2)} \cap W^{\d(K)}} (\dot w \i, 1)(L_{I_2})_{diag} (\dot w
\dot v, B \cap L_{I_2}) \cdot h_{K, \d}=(L_{\d(I_2)}, L_{I_2})
\cdot h_{K, \d}$. Hence $\sqcup_{v \in W_{\d(I_2)} \cap W^{\d(K)}}
(L_{I_2})_{diag} (\dot w \dot v, B \cap L_{I_2}) \cdot h_{K,
\d}=(L_{I_2}, L_{I_2}) \cdot h_{K, \d}$.

In the general case, Consider $\pi: (L_{\d(I_2)}, L_{I_2}) \cdot
h_{K, \d}@>>>\overline{L_{\d(I_2)}/Z(L_{\d(I_2)})} g_0$ defined by
$\pi \bigl( (l_1, l_2) h_{K, \d} \bigr)=(l_1, l_2) \cdot (h_{\d(K)
\cap \d(I_2)} g_0)$ for $l_1 \in L_{\d(I_2)}, l_2 \in L_{I_2}$.
Here $h_{\d(K) \cap \d(I_2)}$ on the right side is the base point
in $\overline{L_{\d(I_2)}/Z(L_{\d(I_2)})}$ that corresponds to
$\d(K) \cap \d(I_2)$. It is easy to see that the morphism is
well-defined. Now define the $T$-action on $(L_{\d(I_2)}, L_{I_2})
\cdot h_{K, \d}$ by $t \cdot \bigl( (l_1, l_2) h_{K, \d} \bigr)=(t
l_1, l_2) h_{K, \d}$ for $t \in T$ and $l_1 \in L_{\d(I_2)}, l_2
\in L_{I_2}$. Then $T$ acts transitively on $\pi \i(a)$ for any $a
\in (L_{\d(I_2)}, L_{I_2}) \cdot (h_{\d(K) \cap \d(I_2)} g_0)$.
Now
$$\eqalignno{& \sqcup_{v \in W_{\d(I_2)} \cap W^{\d(K)}}
\pi \bigl((\dot w \i, 1)(L_{I_2})_{diag} \cdot (\dot w \dot v , B
\cap L_{I_2}) \cdot h_{K, \d} \bigr) \cr &=\sqcup_{v \in
W_{\d(I_2)} \cap W^{\d(K)}} (\dot w \i, 1) (L_{I_2})_{diag} (\dot
w \dot v, B \cap L_{I_2}) \cdot (h_{\d(K) \cap \d(I_2)} g_0) \cr
&=(L_{\d(I_2)}, L_{I_2}) \cdot (h_{\d(K) \cap \d(I_2)} g_0).}$$
Moreover $\sqcup_{v \in W_{\d(I_2)} \cap W^{\d(K)}} (\dot w \i,
1)(L_{I_2})_{diag} \cdot (\dot w \dot v , B \cap L_{I_2}) \cdot
h_{K, \d}$ is stable under $T$-action. Thus $\sqcup_{v \in
W_{\d(I_2)} \cap W^{\d(K)}} (\dot w \i, 1)(L_{I_2})_{diag} \cdot
(\dot w \dot v , B \cap L_{I_2}) \cdot h_{K, \d}=(L_{\d(I_2)},
L_{I_2}) \cdot h_{K, \d}$. The lemma is proved. \qed

\proclaim{Proposition 5.5} Let $(J, w) \in \ci_{\d}$ and $u \in
W_{\d}(J, w)$. Set $I_1=I_1(J, u, \d)$, $I_2=I_2(J, u, \d)$ and
$L^{(J, w, \d)}_u=\sqcup_{K \subset I_1} (L_{I_2}, L_{I_2}) (\dot
u, 1) \cdot h_{K, \d}$. Then we have that

(1) $L^{(J, w, \d)}_u$ is a fibre bundle over
$\overline{L_{I_2}/Z(L_{I_2})}$ with fibres isomorphic to an
affine space of dimension $|I_1|-|I_2|$.

(2) $X^{(J, w, \d)}_u=G_{diag} \cdot L^{(J, w, \d)}_u$ is
isomorphic to $G \times_{P_{I_2}} \bigl((P_{I_2})_{diag} \cdot
L^{(J, w, \d)}_u \bigr)$.

(3) $(P_{I_2})_{diag} \cdot L^{(J, w, \d)}_u=(B \times B) \cdot
L^{(J, w, \d)}_u \cong (U \cap {}^{\dot w^{I_2}_0 \dot u \dot w_0}
U^-) \times L^{(J, w, \d)}_u$.

\endproclaim

Proof. For part (1), note that $L^{(J, u, \d)}_u=\sqcup_{K \subset
I_1} (\dot u, 1) (L_{\d(I_2)}, L_{I_2}) \cdot h_{K, \d}=(\dot u
L_{\d(I_2)}, L_{I_2}) \cdot \bigl(\sqcup_{K \subset I_1} (T, 1)
h_{K, \d} \bigr)$ is a variety. Consider the morphism $$\pi':
\sqcup_{K \subset I_1} (L_{\d(I_2)}, L_{I_2}) \cdot h_{K,
\d}@>>>\overline{L_{\d(I_2)}/Z(L_{\d(I_2)})} g_0$$ defined by
$\pi' \bigl( (l_1, l_2) h_{K, \d} \bigr)=(l_1, l_2) \cdot
(h_{\d(K) \cap \d(I_2)} g_0)$ for $l_1 \in L_{\d(I_2)}, l_2 \in
L_{I_2}$. It is easy to see that $\pi'$ is well defined and is a
locally trivial fibration with fibers isomorphic to an affine
space of dimension $|I_1|-|I_2|$.

Let $v \in W_{\d(I_2)}$. For $K \subset J$, if $\Ad(u v) \d(K)=K$,
then $\Ad(u) \Phi_{\d(K)}=\Ad(u v \i u \i) \Phi_K$. Since $u v \i
u \i \in W_{I_2}$, we have that $\Ad(u) \Phi_{\d(K)} \subset
\Phi_{K \cup I_2}$. Thus $\Ad(u) \Phi_{\d(K \cup I_2)} \subset
\Phi_{K \cup I_2}$. By the maximal property of $I_2$, $K \cup I_2
\subset I_2$. Thus $I_2(J, u v, \d) \subset I_2$. Therefore,

$$\eqalignno{G_{diag} \cdot L^{(J, w, \d)}_u &=G_{diag}
\bigl( \sqcup_{K \subset I_1} (L_{I_2}, L_{I_2}) (\dot u, 1) \cdot
h_{K, \d} \bigr) \cr &=G_{diag} \bigl( \sqcup_{K \subset I_1}
\sqcup_{v \in W_{\d(I_2)} \cap W^{\d(K)}} (L_{I_2})_{diag} (\dot u
\dot v, B \cap L_{I_2(J, u, \d)}) \cdot h_{K, \d} \bigr) \cr
&=\sqcup_{K \subset I_1} \sqcup_{v \in W_{\d(I_2)} \cap W^{\d(K)}}
G_{diag} (\dot u \dot v, B \cap L_{I_2}) \cdot h_{K, \d} \cr
&=\sqcup_{K \subset I_1} \sqcup_{v \in W_{\d(I_2)}} Z^{u v}_{J,
\d}=X^{(J, w, \d)}_u.}$$

Assume that $(g, g) a=b$ for some $g \in G$ and $a, b \in L^{(J,
u, \d)}_u$. Then $a, b$ are in the same $G$ orbit. Note that any
element in $L^{(J, u, \d)}_u$ is conjugate by $L_{I_2}$ to an
element of the form $(\dot u \dot v, l) h_{K, \d}$ with $v \in
W_{\d(I_2)}$, $K \subset I_1(J, u v, \d)$ and $l \in L_{I_2} \cap
B$. Moreover, $(\dot u \dot v, L_{I_2} \cap B) \cdot h_{K, \d}
\subset Z^{u v}_{K, \d}$. Thus if $v_1 \neq v_2$ or $K_1 \neq
K_2$, then for any $l, l' \in L_{I_2} \cap B$, $(\dot u \dot v_1,
l) \cdot h_{K_1, \d}$ and $(\dot u \dot v_2, l') \cdot h_{K_2,
\d}$ are not in the same $G$ orbit. Thus $(g, g) (\dot u \dot v,
l_1) \cdot h_{K, \d}=(\dot u \dot v, l_2) \cdot h_{K, \d}$ for
some $v \in W_{\d(I_2)}$, $K \subset I_1(J, u v, \d)$ and $l_1,
l_2 \in L_{I_2} \cap B$. By 1.12, $g \in P_{I_2(K, u v, \d)}$.
Since $I_2(K, u v, \d) \subset I_2(J, u v, \d) \subset I_2$, we
have that $g \in P_{I_2}$. By 1.9, $X^{(J, w, \d)}_u \cong G
\times_{P_{I_2}} \bigl((P_{I_2})_{diag} \cdot L^{(J, w, \d)}_u
\bigr)$. Part (2) is proved.

For part (3), it is easy to see that $(P_{I_2(J, u, \d)})_{diag}
\cdot L^{(J, w, \d)}_u \subset (B \times B) \cdot L^{(J, w,
\d)}_u$. On the other hand, $$\eqalignno{(B \times B) \cdot L^{(J,
w, \d)}_u &=(U_{P_{I_2}}, U_{P_{I_2}}) (L_{I_2})_{diag}
\bigl(\sqcup_{v \in W_{\d(I_2)}} \sqcup_{K \subset I_2(J, u v,
\d)} (\dot u \dot v, B) \cdot h_{K, \d} \bigr) \cr
&=(L_{I_2})_{diag} (U_{P_{I_2}}, U_{P_{I_2}}) \bigl(\sqcup_{v \in
W_{\d(I_2)}} \sqcup_{K \subset I_2(J, u v, \d)} (\dot u \dot v, B)
\cdot h_{K, \d} \bigr).}$$

By 1.12, $(U_{P_{I_2}}, U_{P_{I_2}}) (\dot u \dot v, B) \cdot
h_{K, \d}=(B \times B) (\dot u \dot v, 1) \cdot h_{K, \d} \subset
(P_{I_2(K, u v, \d)})_{diag} \cdot (L_{I_2(K, u v, \d)}, L_{I_2(K,
u v, \d)}) (\dot u \dot v, 1) h_{K, \d}$. We have showed that
$I_2(K, u v, \d) \subset I_2$. Hence $(U_{P_{I_2}}, U_{P_{I_2}})
(\dot u \dot v, B) \cdot h_{K, \d} \subset (P_{I_2(J, u,
\d)})_{diag} \cdot L^{(J, w, \d)}_u$. Therefore, $(P_{I_2})_{diag}
\cdot L^{(J, w, \d)}_u=(B \times B) \cdot L^{(J, w, \d)}_u$.

Consider the morphism $\pi: (U \cap {}^{\dot w^{I_2}_0 \dot u \dot
w_0} U^-) \times L^{(J, w, \d)}_u@>>>(B \times B) \cdot L^{(J, w,
\d)}_u$ defined by $\pi(b, l)=(b, 1) \cdot l$ for $b \in U \cap
{}^{\dot w^{I_2}_0 \dot u \dot w_0} U^-$ and $l \in L^{(J, w,
\d)}_u$. By the similar argument as we did in 1.10, we can show
that $\pi$ is an isomorphism. \qed

\proclaim{Corollary 5.6} We keep the notation of 5.5. If moreover,
$I_2=\varnothing$, then $X^{(J, w, \d)}_u$ admits a cellular
decomposition.
\endproclaim

Proof. If $I_2=\varnothing$, then $L^{(J, w, \d)}_u$ is an affine
space. Thus $X^{(J, w, \d)}_u$ is isomorphic to $G \times B C$,
where $C=(B \times B) \cdot L^{(J, w, \d)}_u$. By part (3) of 5.5,
$C$ is an affine space. It is easy to see that $B$ acts linearly
on $C$. Therefore $X^{(J, w, \d)}_u$ is a vector bundle over
$\cb$. Note that $\cb$ admits a cellular decomposition. By a
well-known result (see \cite{Q}, \cite{Su} or \cite{VS}), $X^{(J,
w, \d)}_u$ admits a cellular decomposition. \qed

\subhead 5.7\endsubhead For $w_1, w_2 \in W_{\d}(J, w)$, we say
$w_2 \le' w_1$ if there exists $w_1=x_0, x_1, \cdots, x_n=w_2$,
$v_i \in \d(I_2(J, x_{i+1}, \d))$ for all $i$, such that $x_{i+1}
v_i \ge_{I_1(J, x_i, \d), \d} x_i$.

By 4.5, $\overline{X^{(J, w, \d)}_{u_1}} \cap X^{(J, w,
\d)}_{u_2}=\varnothing$ if $u_2 \nleqslant' u_1$. hence if $\le'$
is a partial order on $W_{\d}(J, w)$, then $\overline{Z^w_{J,
\d}}=\sqcup_{u \in W_{\d}(J, w)} X^{(J, w, \d)}_u$ is an
$\a$-partition. We will show that $\le'$ is a partial order if
$\overline{Z^w_{J, \d}}$ contains finitely many $G$-orbits.

\proclaim{Lemma 5.8} Let $J \subset I$, $u \in W$, $w \in W^J$ and
$v \in W_J$. Assume that $u w v=w' v'$ for some $w' \in W^J$ and
$v' \in W_J$. If $l(u w v)=l(w v)-l(u)$, then $w' \le w$. If
moreover, $w'=w$, then $\Ad(w \i) \supp(u) \subset J$.
\endproclaim

Proof. If $u=s_i$ for some $i \in J$ and $l(s_j w v)=l(w v)-1$,
then either $s_i w<w$ and $s_i w \in W^J$ or $s_i w=w s_j$ for
some $j \in J$. It is easy to check that the statement holds in
both cases.

The general case can be proved by induction on $l(u)$. \qed

\proclaim{Lemma 5.9} If $w_1, w_2 \in W_{\d}(J, w)$ with $w_1 \le'
w_2$ and $w_2 \le' w_1$ and $I_2(J, w_1, \d)=I_2(J, w_2,
\d)=\varnothing$, then $w_1=w_2$.
\endproclaim

Proof. We will prove the case: if $w_1 \ge_{I_1(J, w_2, \d), \d}
w_2$, $w_2 \ge_{I_1(J, w_1, \d), \d} w_1$ and $I_2(J, w_1,
\d)=\varnothing$, then $w_1=w_2$. The general case can be proved
in the similar way.

We argue by induction on $|J|$. Since $l(w_1) \ge l(w_2)$ and
$l(w_2) \ge l(w_1)$, we have that $l(w_1)=l(w_2)$. Thus $w_1=u_2
\i w_2 \d(u_2)$ and $w_2=u_1 \i w_1 \d(u_1)$ for some $u_1 \in
W_{I_1(J, w_1, \d)}$ and $u_2 \in W_{I_1(J, w_2, \d)}$. By
induction hypothesis, it suffices to prove the case when
$J=\supp(u_1) \cup \supp(u_2)$.

We have that $w_1=w'_1 \d(v_1)$ and $w_2=w'_2 \d(v_2)$ for some
$w'_1, w'_2 \in W^{\d(J)}$ and $v_1, v_2 \in W_J$. Note that $w'_1
\d(v_1)=u_2 \i w'_2 \d(v_2 u_2)$ and $l(u_2 \i w'_2 \d(v_2
u_2))=l(w'_2 \d(v_2 u_2))-l(u_2)$. By 5.8, $w'_1 \le w'_2$.
Similarly $w'_2 \le w'_1$. Therefore $w'_1=w'_2$. By 5.8,
$\Ad(w'_2) \i \supp(u_2) \subset \d(J)$ and $\Ad(w'_1) \i
\supp(u_1) \subset \d(J)$. Therefore $\Ad(w'_1) \i J \subset
\d(J)$. Hence $\Ad(w_1) \i \Phi_J=\Phi_{\d(J)}$. Since $I_2(J,
w_1, \d)=\varnothing$, we have that $J=\varnothing$. Therefore
$w_1 \ge w_2$ and $w_2 \ge w_1$. Thus $w_1=w_2$. The case is
proved. \qed

${}$

As a summary, we have the following result.

\proclaim{Theorem 5.10} If $\overline{Z^w_{J, \d}}$ contains only
finitely many $G$-orbits, then it has a cellular decomposition.
\endproclaim

Proof. If $\overline{Z^w_{J, \d}}$ contains only finitely many
$G$-orbits, then $I_2(J, u, \d)=\varnothing$ for $u \ge_{J, \d}
w$. In this case, $W_{\d}(J, w)=\{u \in W \mid u \ge_{J, \d} w\}$.
We have that $$\overline{Z^w_{J, \d}}=\sqcup_{u \in W_{\d}(J, \d)}
X^{(J, w, \d)}_u.$$

By 5.7 and 5.9, the partition is an $\a$-partition. Thus by 5.6,
$\overline{Z^w_{J, \d}}$ has a cellular decomposition. \qed

\head Acknowledgements. \endhead We thank George Lusztig for
suggesting the problem and for many helpful discussions. We thank
T. A. Springer and David Vogan for some useful comments. We also
thank Jiang-hua Lu for point out a mistake in the previous version
and suggested the reference \cite{SL}.


\Refs

\widestnumber\key{DP}

\ref\key{DP} \by C. De Concini and C.Procesi\paper Complete
symmetric varieties\inbook Invariant theory (Montecatini 1982),
Lect. Notes Math.\vol 996\pages 1-44\publ Springer \yr 1983
\endref

\ref\key{H} \by X. He \paper Unipotent variety in the group
compactification \jour Adv. in Math., in press\endref

\ref\key{L1} \by G. Lusztig\paper Total positivity in reductive
groups \jour Lie Theory and Geometry: in honor of Bertram Kostant,
Progress in Math. \vol 123\pages 531-568\publ
Birkh\"auser\publaddr Boston \yr 1994\endref

\ref\key{L2} \by G. Lusztig\book Hecke algebras with unequal
parameters, CRM Monograph Series, 18 \publ American Mathematical
Society \publaddr Providence, RI \yr 2003\endref

\ref\key{L3} \by G. Lusztig\paper Parabolic character sheaves I
\jour Moscow Math.J \vol 4 \yr 2004 \pages 153-179
\endref

\ref\key{L4} \by G. Lusztig\paper Parabolic character sheaves II
\jour Moscow Math.J \vol 4 \yr 2004 \pages 869-896
\endref

\ref\key{Q} \by D. Quillen\paper Projective modules over
polynomial rings \jour Invent. Math. \vol 36 \yr 1976 \pages
167-171 \endref

\ref\key{S} \by T. A. Springer \paper Intersection cohomology of
$B \times B$-orbits closures in group compactifications, \jour J.
Alg. \vol 258 \yr 2002 \pages 71-111
\endref

\ref\key{SL} \by P. Slodowy \book Simple singularities and simple
algebraic groups, Lecture Notes in Mathematics, 815 \publ Springer
\publaddr Berlin \yr 1980 \endref


\ref\key{Su} \by A. A. Suslin \paper Projective modules over
polynomial rings are free (Russian) \jour Dokl. Akad. Nauk SSSR
\vol 229 \yr 1976 \pages no. 5, 1063-1066 \endref

\ref\key{VS} \by L. N. Vaserstein and A. A. Suslin \paper Serre's
problem on projective modules over polynomial rings, and algebraic
$K$-theory (Russian) \jour Izv. Akad. Nauk SSSR Ser. Mat. \vol 40
\yr 1976 \pages no. 5, 993--1054, 1199
\endref

\endRefs
\enddocument